\documentclass[reqno, a4paper]{amsart}
\usepackage{amsmath,amssymb,dsfont,verbatim}

\usepackage[raggedright]{titlesec}

\titleformat{\chapter}[display]
{\normalfont\huge\bfseries}{\chaptertitlename\\thechapter}{20pt}{\Huge}
\titleformat{\section}
{\normalfont\Large\bfseries\center}{\thesection}{1em}{}
\titleformat{\subsection}
{\normalfont\large\bfseries}{\thesubsection}{1em}{}
\titleformat{\subsubsection}[runin]
{\normalfont\normalsize\bfseries}{\thesubsubsection}{1em}{}
\titleformat{\paragraph}[runin]
{\normalfont\normalsize\bfseries}{\theparagraph}{1em}{}
\titleformat{\subparagraph}[runin]
{\normalfont\normalsize\bfseries}{\thesubparagraph}{1em}{}
\titlespacing*{\chapter} {0pt}{50pt}{40pt}
\titlespacing*{\section} {0pt}{3.5ex plus 1ex minus .2ex}{2.3ex plus .2ex}
\titlespacing*{\subsection} {0pt}{3.25ex plus 1ex minus .2ex}{1.5ex plus .2ex}
\titlespacing*{\subsubsection}{0pt}{3.25ex plus 1ex minus .2ex}{1.5ex plus .2ex}
\titlespacing*{\paragraph} {0pt}{3.25ex plus 1ex minus .2ex}{1em}
\titlespacing*{\subparagraph} {\parindent}{3.25ex plus 1ex minus .2ex}{1em}

\input xypic
\xyoption{all}



\newtheorem{theorem}{Theorem}[section]
\newtheorem{lemma}[theorem]{Lemma}
\newtheorem{proposition}[theorem]{Proposition}
\newtheorem{corollary}[theorem]{Corollary}

\theoremstyle{definition}
\newtheorem{definition}[theorem]{Definition}

\newtheorem{notation}[theorem]{Notation}
\newtheorem{example}[theorem]{Example}

\theoremstyle{remark}
\newtheorem{remark}[theorem]{Remark}

\DeclareMathOperator{\Aut}{Aut}
\DeclareMathOperator{\ide}{id}

\DeclareMathOperator{\HH}{HH}
\DeclareMathOperator{\Ho}{H}
\DeclareMathOperator{\Hom}{Hom}
\DeclareMathOperator{\sg}{sg}

\newcommand{\ov}{\overline}
\newcommand{\ot}{\otimes}

\newcommand{\wt}{\widetilde}
\newcommand{\ep}{\epsilon}
\newcommand{\De}{\Delta}

\def\xcirc{\objectmargin{0.1pc}\def\objectstyle{\sssize}\diagram \squarify<1pt>{}\circled \enddiagram}

\begin{document}

\title{Universal deformation formulas and braided module algebras}

\author{Jorge A. Guccione}
\address{Departamento de Matem\'atica\\ Facultad de Ciencias Exactas y Naturales, Pabell\'on~1, Ciudad Universitaria\\ (1428) Buenos Aires, Argentina.}
\email{vander@dm.uba.ar}
\thanks{Supported by  UBACYT 095, PIP 112-200801-00900 (CONICET) and PUCP-DAI-2009-0042}

\author{Juan J. Guccione}
\address{Departamento de Matem\'atica\\ Facultad de Ciencias Exactas y Naturales\\ Pabell\'on~1, Ciudad Universitaria\\ (1428) Buenos Aires, Argentina.}
\email{jjgucci@dm.uba.ar}
\thanks{Supported by  UBACYT 095 and PIP 112-200801-00900 (CONICET)}

\author{Christian Valqui}
\address{Pontificia Universidad Cat\'olica del Per\'u - Instituto de Matem\'atica y Ciencias Afi\-nes, Secci\'on Matem\'aticas, PUCP, Av. Universitaria 1801, San Miguel, Lima 32, Per\'u.}
\email{cvalqui@pucp.edu.pe}
\thanks{Supported by PUCP-DAI-2009-0042,  Lucet 90-DAI-L005, SFB 478 U. M\"unster, Konrad Adenauer Stiftung.}

\thanks{The second author thanks the appointment as a visiting professor ``C\'atedra Jos\'e Tola Pasquel'' and the hospitality during his stay at the PUCP}

\begin{abstract} We study formal deformations of a crossed product $S(V)\#_f G$, of a polynomial algebra with a group, induced from a universal deformation formula introduced by Witherspoon.  These  deformations arise from braided actions of Hopf algebras generated by automorphisms and skew derivations. We show that they are nontrivial in the characteristic free context, even if $G$ is infinite, by showing that their infinitesimals are not coboundaries. For this we construct a new complex which computes the Hochschild cohomology of $S(V)\#_f G$.
\end{abstract}

\maketitle

\section*{Introduction} In~\cite{G-Z} Giaquinto and Zhang develop the notion of a universal deformation formula based on a bialgebra $H$, extending earlier formulas based on universal enveloping algebras of Lie algebras. Each one of these formulas is called universal because it provides a formal deformation for any $H$-module algebra. In the same paper the authors construct the first family of such formulas based on noncommutative bialgebras, namely the enveloping algebras of central extensions of a Heisenberg Lie algebra $L$. Another of these formulas, based on a Hopf algebra $H_q$ over~$\mathds{C}$, where $q\in~\mathds{C}^{\times}$ is a parameter, generated by group like elements $\sigma^{\pm 1}$ and two skew primitive elements $D_1,D_2$, were obtained in the generic case by the same authors, but were not published. In~\cite{W} the author generalizes this formula to include the case where $q$ is a root of unity, and she uses it to construct formal deformations of a crossed product $S(V)\#_f G$, where $S(V)$ is the polynomial algebra and the group $G$ acts linearly on $V$. More precisely, she deals with deformations whose infinitesimal sends $V\ot V$ to $S(V)w_g$, where $g$ is a central element of $G$.

In this paper we prove that some results established in~\cite{W} under the hypothesis that $G$ is a finite group, remain valid for arbitrary groups, and with $\mathds{C}$ replaced by an arbitrary field. For instance we show that the determinant of the action of $g$ on $V$ is always~$1$. Moreover, we do not only consider standard $H_q$-module algebra structures on $S(V)\#_f G$, but also the more general ones introduced in~\cite{G-G1}, and we work with actions which depend on two central elements $g_1$ and $g_2$ of $G$ and two polynomials $P_1$ and $P_2$. When the actions are the standard ones, $g_1 = 1$ and $P_1 = 1$, we obtain the case considered in~\cite{W}. Finally, in Subsection \ref{Second case} we show how to extend the explicit formulas obtained previously, to non central $g_1$ and $g_2$. As was noted by Witherspoon, these formulas necessarily involve all components of $S(V)\#_f G$ corresponding to the elements of a union of conjugacy classes of $G$.

\smallskip

The paper is organized as follows: in the first section we review the concept of braided module algebra introduced in~\cite{G-G1}, we adapt the notion of Universal Deformation Formula (UDF) to the braided context, and we show that each one of these formulas produces a deformation on any braided $H$-module algebra whose transposition (see Definition~\ref{def: transp alg}) satisfy a suitable hypothesis. We remark that, when the bialgebra $H$ is standard, the use of braided module algebra gives rise to more deformations than the ones obtained using only module algebras, because the transposition can be different from the flip. With this in mind, although we are going to work with the standard Hopf algebra $H_q$, we establish the basic properties of UDF's in the braided case, because it is the most appropriate setting to deal with arbitrary transpositions. In the second section we recall the definitions of the Hopf algebra $H_q$ and of the UDF $\exp_q$ considered in~\cite[Section~3]{W}, which we are going to study. We also introduce the concept of a good transposition of $H_q$ on an algebra $A$, and we study some of its properties. Perhaps the most important result in this section is Theorem~\ref{algunas acciones de H_q sobre algebras}, in which we obtain a description of all the $H_q$-module algebras $(A,s)$, with $s$ a good transposition. This is the first of several results in which we give a systematic account of the necessary and sufficient conditions that an algebra (in general a crossed product $S(V)\#_f G$) must satisfy in order to support a braided $H_q$-module algebra structure satisfying suitable hypothesis. In Section~4 of~\cite{W}, using the UDF $\exp_q$ the author constructs a large family of deformations whose infinitesimal sends $V\ot V$ to $S(V)w_g$, where $g$ is a central element of $G$. Using cohomological methods she proves that if $G$ is finite, these deformations are non trivial, that the action of $g$ on $V$ has determinant~$1$ and that the codimension of ${}^gV$ is~$0$ or~$2$. In the first part of Section~3 we study a larger family of deformations and we prove that the last two results hold for this family even if $G$ is infinite and the characteristic of $k$ is non zero. Finally, in Section~4 we show that, under very general hypothesis, the deformations constructed in the previous section are non trivial. Once again, we do not assume characteristic zero, nor that the group $G$ is finite. One of the interesting points in this paper is the method developed to deal with the cohomology of $S(V)\#_f G$ when $k[G]$ is non semisimple.  As far as we know it is the first time that this type of cochain complexes is used to prove the non triviality of a Hochschild cocycle.

\section{Preliminaries}
After introducing some basic notations we recall briefly the concepts of braided bialgebra and braided Hopf algebra following the presentation given in~\cite{T1} (see also~\cite{T2}, \cite{Ly1}, \cite{F-M-S}, \cite{A-S}, \cite{D}, \cite{So} and \cite{B-K-L-T}). Then we review the notion of braided module algebra introduced in~\cite{G-G1}, we recall the concept of Universal Deformation Formula based on a bialgebra $H$, due to Giaquinto an Zhang, and we show that such a UDF produces a formal deformation when it is applied to an $H$-braided module algebra, satisfying suitable hypothesis, generalizing slightly a result in~\cite{G-Z}.

\smallskip

In this paper $k$ is a field, $k^{\times} = k\setminus\{0\}$, all the vector spaces are over $k$, and $\ot=\ot_k$. Moreover we will use the usual notation $(i)_q = 1 + q +\cdots + q^{i-1}$ and $(i)!_q = (1)_q\cdots (i)_q$, for $q\in k^{\times}$ and $i\in\mathds{N}$.

\smallskip

Let $V$, $W$ be vector spaces and let $c\colon V\ot W\to W\ot V$ be a $k$-linear map. Recall that:

\begin{itemize}

\item[-] If $V$ is an algebra, then $c$ is compatible with the algebra structure of $V$ if $c \xcirc (\eta\ot W)= W\ot\eta$ and $c\xcirc(\mu\ot W) = (W\ot\mu)\xcirc(c\ot V)\xcirc(V\ot c)$, where $\eta\colon k\to V$ and $\mu\colon V\ot V\to V$ denotes the unit and the multiplication map of $V$, respectively.

\smallskip

\item[-] If $V$ is a coalgebra, then $c$ is compatible with the coalgebra structure of $V$ if $(W\ot \ep)\xcirc c = \ep\ot W$ and $(W\ot\De)\xcirc c = (c\ot V)\xcirc (V\ot c)\xcirc (\De\ot W)$, where $\epsilon\colon V\to k$ and $\Delta\colon V\to V\ot V$ denotes the counit and the comultiplication map of $V$, respectively.

\end{itemize}

Of course, there are similar compatibilities when $W$ is an algebra or a coalgebra.

\subsection{Braided bialgebras and braided Hopf algebras}

\begin{definition} A {\em braided bialgebra} is a vector space $H$ endowed with an algebra structure, a coalgebra structure and a braiding operator $c\in\Aut_k(H^{\ot 2})$ (called the {\em braid} of $H$), such that $c$ is compatible with the algebra and coalgebra structures of $H$, $\De\xcirc \mu = (\mu\ot\mu)\xcirc(H\ot c\ot H)\xcirc(\De\ot\De)$, $\eta$ is a coalgebra morphism and $\ep$ is an algebra morphism. Furthermore, if there exists a $k$-linear map $S\colon H\to H$, which is the inverse of the identity map for the convolution product, then we say that $H$ is a {\em braided Hopf algebra} and we call $S$ the {\em antipode} of $H$.\label{de1.1}
\end{definition}

Usually $H$ denotes a braided bialgebra, understanding the structure maps, and $c$ denotes its braid. If necessary, we will use notations as $c_H$, $\mu_H$, etcetera.

\begin{remark} Assume that $H$ is an algebra and a coalgebra and $c\in\Aut_k(H^{\ot 2})$ is a solution of the braiding equation, which is compatible with the algebra and coalgebra structures of $H$. Let $H\ot_c H$ be the algebra with underlying vector space $H^{\ot 2}$ and multiplication map given by $\mu_{H\ot_c H}:= (\mu\ot\mu)\xcirc (H\ot c\ot H)$. It is easy to see that $H$ is a braided bialgebra with braid $c$ if and only if $\De\colon H\to H\ot_c H$ and $\ep\colon H\to k$ are morphisms of algebras.\label{re1.2}
\end{remark}

\begin{definition} Let $H$ and $L$ be braided bialgebras. A map $g\colon H\to L$ is a {\em morphism of braided bialgebras} if it is an algebra homomorphism, a coalgebra homomorphism and $c\xcirc(g\ot g)=(g\ot g)\xcirc c$.\label{de1.3}
\end{definition}

Let $H$ and $L$ be braided Hopf algebras. It is well known that if $g\colon H\to L$ is a morphism of braided bialgebras, then $g\xcirc S = S\xcirc g$.
\subsection{Braided module algebras}
\begin{definition}\label{def: H-braided space} Let $H$ be a braided bialgebra. A {\em left $H$-braided space} $(V,s)$ is a vector space $V$, endowed with a bijective $k$-linear map $s\colon H\ot V\to V\ot H$, which is compatible with the bialgebra structure of $H$ and satisfies
$$
(s\ot H)\xcirc (H\ot s)\xcirc (c\ot V) = (V\ot c)\xcirc (s\ot H)\xcirc (H\ot s)
$$
(compatibility of $s$ with the braid). Let $(V',s')$ be another left $H$-braided space. A $k$-linear map $f\colon V\to V'$ is said to be a {\em morphism of left $H$-braided spaces}, from $(V,s)$ to $(V',s')$, if $(f\ot H)\xcirc s = s'\xcirc (H\ot f)$.
\end{definition}

We let $\mathcal{LB}_H$ denote the category of all left $H$-braided spaces. It is easy to check that this is a monoidal category with:

\smallskip

\begin{itemize}

\item[-] unit $(k,\tau)$, where $\tau\colon H\ot k\to k\ot H$ is the flip,

\smallskip

\item[-] tensor product $(V,s_V)\ot(U,s_U) := (V\ot U, s_{V\ot U})$, where $s_{V\ot U}$ is the map $s_{V\ot U}:= (V\ot s_U) \xcirc(s_V\ot U)$,

\smallskip

\item[-] the usual associativity and unit constraints.

\end{itemize}

\begin{definition}\label{def: braided alg} We will say that $(A,s)$ is a {\em left $H$-braided algebra} or simply a {\em left $H$-algebra} if it is an algebra in $\mathcal{LB}_H$.
\end{definition}

We let $\mathcal{ALB}_H$ denote the category of left $H$-braided algebras.

\begin{definition}\label{def: transp alg} Let $A$ be an algebra. A {\em left transposition} of $H$ on $A$ is a bijective map $s\colon H\ot A\to A\ot H$, satisfying:

\smallskip

\begin{enumerate}

\item $(A,s)$ is a left $H$-braided space,

\smallskip

\item $s$ is compatible with the algebra structure of $A$.

\end{enumerate}

\end{definition}

\begin{remark}\label{re: alg in LB_H} A left $H$-braided algebra is a pair $(A,s)$ consisting of an algebra $A$ and a left transposition $s$ of $H$ on $A$. Let $(A',s')$ be another left $H$-braided algebra. A map $f\colon A\to A'$ is a morphism of left $H$-braided algebras, from $(A,s)$ to $(A',s')$, if and only if it is a morphism of standard algebras and $(f\ot H)\xcirc s = s'\xcirc (H\ot f)$.
\end{remark}

Note that $(H,c)$ is an algebra in $\mathcal{LB}_H$. Hence, one can consider left and right $(H,c)$-modules in this monoidal category.

\begin{definition}\label{def: H-braided module} We will say that $(V,s)$ is a {\em left $H$-braided module} or simply a {\em left $H$-module} to mean that it is a left $(H,c)$-module in $\mathcal{LB}_H$.
\end{definition}

We let ${}_H(\mathcal{LB}_H)$ denote the category of left $H$-braided modules.

\begin{remark}\label{re: H-braided module} A left $H$-braided space $(V,s)$ is a left $H$-module if and only if $V$ is a standard left $H$-module and
$$
s\xcirc (H\ot\rho) = (\rho\ot H)\xcirc (H\ot s)\xcirc (c\ot V),
$$
where $\rho$ denotes the action of $H$ on $V$. Furthermore, a map $f\colon V\to V'$ is a {\em morphism of left $H$-modules}, from $(V,s)$ to $(V',s')$, if and only if it is $H$-linear and $(f\ot H)\xcirc s = s'\xcirc (H\ot f)$.
\end{remark}

Given left $H$-modules $(V,s_V)$ and $(U,s_U)$, with actions $\rho_V$ and $\rho_U$ respectively, we let $\rho_{V\ot U}$ denote the diagonal action
$$
\rho_{V\ot U}:= (\rho_V\ot\rho_U)\xcirc(H\ot s_V\ot U)\xcirc (\De_H\ot V\ot U).
$$
The following proposition says in particular that $(k,\tau)$ is a left $H$-module via the trivial action and that $(V,s_V)\ot (U,s_U)$ is a left $H$-module via $\rho_{V\ot U}$.

\begin{proposition}[see \cite{G-G1}]\label{prop: _H sub LB H es monoidal} The category ${}_H(\mathcal{LB}_H)$, of left $H$-braided modules, endowed with the usual associativity and unit constraints, is monoidal.
\end{proposition}

\begin{definition}\label{def: H-braided mod alg} We say that $(A,s)$ is a {\em left $H$-braided module algebra} or simply a {\em left $H$-module algebra} if it is an algebra in ${}_H(\mathcal{LB}_H)$.
\end{definition}

We let ${}_H(\mathcal{ALB}_H)$ denote the category of left $H$-braided module algebras.

\begin{remark}\label{re: car de H-braid mod alg} $(A,s)$ is a left $H$-module algebra if and only if the following facts hold:

\begin{enumerate}

\smallskip

\item $A$ is an algebra,

\smallskip

\item $s$ is a left transposition of $H$ on $A$,

\smallskip

\item $A$ is a standard left $H$-module,

\smallskip

\item $s\xcirc (H\ot\rho) = (\rho\ot H)\xcirc (H\ot s)\xcirc (c\ot A)$,

\smallskip

\item $\mu_A\xcirc (\rho\ot\rho)\xcirc (H\ot s\ot A)\xcirc (\De_H\ot A\ot A) = \rho \xcirc (H\ot\mu_A)$,

\smallskip

\item $h\cdot 1 = \ep(h)1$ for all $h\in H$,

\smallskip

\end{enumerate}
where $\rho$ denotes the action of $H$ on $A$. So, $(A,s)$ is a left $H$-module algebra if and only if it is a left $H$-algebra, a left $H$-module and satisfies conditions~(5) and~(6).
\end{remark}

In the sequel, given a map $\rho \colon H\ot A\to A$, sometimes we will write $h\cdot a$ to denotes $\rho(h\ot a)$.

\begin{remark}\label{condiciones para que (A,s) sea un H-módulo álgebra a izquierda} If $X$ generates $H$ as a $k$-algebra, then conditions~(4), (5) and~(6) of the above remark are satisfied if and only if
\begin{align*}
&s(h\ot l\!\cdot\!a) = (\rho\ot H)\xcirc (H\ot s)\xcirc (c\ot A)(h\ot l\ot a),\\
& h\cdot(ab) = \mu_A\xcirc(\rho\ot\rho)\xcirc(H\ot s\ot A)(\Delta(h)\ot a\ot b),\\
& h\!\cdot\!1 = \epsilon (h),
\end{align*}
for all $a,b\in A$ and $h,l\in X$.
\end{remark}

Let $(A',s')$ be another left $H$-module algebra. A map $f\colon A\to A'$ is a {\em morphism of left $H$-module algebras}, from $(A,s)$ to $(A',s')$, if and only if it is an $H$-linear morphism of standard algebras that satisfies $(f\ot H)\xcirc s= s'\xcirc (H\ot f)$.

\subsection{Bialgebra actions and universal deformation formulas}
Most of the results of~\cite[Section~1]{G-Z} remain valid in our more general context, with the same arguments and minimal changes. In particular Theorem~\ref{uso de twisting elements para torcer algebras} below holds.

Let $H$ be a braided bialgebra. Given a left $H$-module algebra $(A,s)$ and an element $F\in H\ot H$, we let $F_l\colon A\ot A\to A\ot A$ denote the map defined by
$$
F_l(a\ot b) := (\rho\ot\rho)\xcirc (H\ot s\ot A)(F\ot a\ot b),
$$
where $\rho\colon H\ot A\to A$ is the action of $H$ on $A$. We let $A_F$ denote $A$ endowed with the multiplication map $\mu_A\xcirc F_l$.

\begin{definition}\label{twisting element} We say that  $F\in H\ot H$ is a {\em twisting element (based on $H$)} if

\begin{enumerate}

\smallskip

\item $(\epsilon\ot\ide)(F) = (\ide\ot\epsilon)(F) = 1$,

\smallskip

\item $[(\Delta\ot\ide)(F)](F\ot 1) = [(\ide\ot\Delta)(F)](1\ot F)$ in $H\ot_c H\ot_c H$,

\smallskip

\item $(c\ot H)\xcirc(H\ot c)(F\ot h) = h\ot F$, for all $h\in H$.

\end{enumerate}

\end{definition}

\begin{theorem}\label{uso de twisting elements para torcer algebras} Let $(A,s)$ be a left $H$-module algebra. If $F\in H\ot H$ is a twisting element such that $(s\ot H)\xcirc(H\ot s)(F\ot a) = a\ot F$, for all $a\in A$, then $A_F$ is an associative algebra with unit $1_A$.
\end{theorem}

The notions of braided bialgebra, left $H$-braided module algebra and twisting element make sense in arbitrary monoidal categories. Here we consider the monoi\-dal category $\mathcal{M}[[t]]$ defined as follows:

\begin{itemize}

\smallskip

\item[-] the objects are the $k[[t]]$-modules of the form $M[[t]]$ where $M$ is a $k$-vector space,

\smallskip

\item[-] the arrows are the $k[[t]]$-linear maps,

\smallskip

\item[-] the tensor product is the completation
$$
M[[t]]\widehat{\ot}_{k[[t]]} N[[t]]
$$
of the algebraic tensor product $M[[t]]\ot_{k[[t]]} N[[t]]$ with respect to the $t$-adic topology,

\smallskip

\item[-] the unit and the associativity constrains are the evident ones.

\smallskip

\end{itemize}

We identify $M[[t]]\widehat{\ot}_{k[[t]]} N[[t]]$ with $(M\ot N)[[t]]$ by the map
$$
\Theta\colon M[[t]]\widehat{\ot}_{k[[t]]} N[[t]]\to (M\ot N)[[t]]
$$
given by $\Theta(mt^i\ot nt^j) := (m\ot n)t^{i+j}$.

\medskip

If $A$ is a $k$-algebra, then $A[[t]]$ is an algebra in $\mathcal{M}[[t]]$ via the multiplication map
$$
\xymatrix @R=-2pt {(A\ot A)[[t]]\rto^-{\mu} & A[[t]]\\
\sum (a_i\ot b_i) t^i \ar@{|->}[0,1] & \sum a_ib_it^i},
$$
where $a_ib_i = \mu_A(a_i\ot b_i)$. The unit map is the canonical inclusion $k[[t]]\hookrightarrow A[[t]]$.

\smallskip

If $H$ is a braided bialgebra over $k$, then $H[[t]]$ is a braided bialgebra in $\mathcal{M}[[t]]$. The multiplication and unit maps are as above. The comultiplication and counits are the maps
$$
\xymatrix @R=-2pt {H[[t]]\rto^-{\Delta} & (H\ot H) [[t]]\\
\sum h_i t^i \ar@{|->}[0,1] & \sum \Delta_H(h_i) t^i}\qquad\text{and}\qquad
\xymatrix @R=-2pt {H[[t]]\rto^-{\epsilon} & k [[t]]\\
\sum h_i t^i \ar@{|->}[0,1] & \sum \epsilon_H(h_i) t^i},
$$
and the braid operator is the map
$$
\xymatrix @R=-2pt {(H\ot H)[[t]]\rto^-{c[[t]]} & (H\ot H) [[t]]\\
\sum (h_i\ot l_i) t^i \ar@{|->}[0,1] & \sum c_H(h_i\ot l_i) t^i}.
$$

\smallskip

If $(A,s)$ is a $H$-module algebra, then $(A[[t]],s[[t]])$, where $s[[t]]$ is the map
$$
\xymatrix @R=-2pt {(H\ot A)[[t]]\rto^-{s[[t]]} & (A\ot H) [[t]]\\
\sum (h_i\ot a_i) t^i \ar@{|->}[0,1] & \sum s(h_i\ot a_i) t^i},
$$
is an $H[[t]]$-module algebra, via
$$
\xymatrix @R=-2pt {(H\ot A)[[t]]\rto^-{\rho} & A [[t]]\\
\sum (h_i\ot a_i) t^i \ar@{|->}[0,1] & \sum \rho_A(h_i\ot a_i) t^i}.
$$

\smallskip

A twisting element based on $H[[t]]$ in $\mathcal{M}[[t]]$ is an element $F\in H[[t]] \widehat{\ot}_{k[[t]]} H[[t]]$ satisfying conditions (1)--(3) of Definition~\ref{twisting element}. It is easy to check that a power series $F=\sum F_i\, t^i\in (H\ot H)[[t]]$ corresponds via $\Theta^{-1}$ to a twisting element if and only if

\begin{enumerate}

\smallskip

\item $(\epsilon \ot \ide)(F_0)=(\ide\ot \epsilon)(F_0)=1$ and $(\epsilon \ot \ide)(F_i) = (\ide\ot \epsilon)(F_i)=0$ for $i\ge 1$,

\smallskip

\item For all $n\ge 0$,
$$
\qquad\qquad\sum_{i+j=n} (\Delta\ot \ide)(F_i)(F_j\ot 1) = \sum_{i+j=n} (\ide\ot \Delta)(F_i) (1\ot F_j) \quad\text{in $H\ot_c H\ot_c H$,}
$$

\smallskip

\item $(c\ot H)\xcirc(H\ot c)(F_n\ot h) = h\ot F_n$, for all $h\in H$ and $n\ge 0$.

\smallskip

\end{enumerate}

We will say that $F$ is an {\em universal deformation formula} (UDF) {\em based on $H$} if, moreover, $F_0=1\ot 1$.

\smallskip

\begin{theorem}\label{uso de UDF para construir deformaciones} Let $(A,s)$ be a left
$H$-module algebra. If $F = \sum F_i\, t^i$ is an UDF based on $H$, such that
$$
(s\ot H)\xcirc(H\ot s)(F_i\ot a) = a\ot F_i\quad\text{for all $i\ge 0$ and $a\in A$,}
$$
then, the construction considered in Theorem~\ref{uso de twisting elements para torcer  algebras}, applied to the left $H[[t]]$-module algebra $(A[[t]],s[[t]])$ introduced above, produces a formal deformation of $A$.
\end{theorem}

\begin{proof} It is immediate.
\end{proof}

\section{$\mathbf{H}_\mathbf{q}$-module algebra structures and deformations}
In this section, we briefly review the construction of the Hopf algebra $H_q$ and the UDF $\exp_q$ based on $H_q$ considered in~\cite{W}, we introduce the notion of a good transposition of $H_q$ on an algebra $A$, and we describe all the braided $H_q$-module algebras whose transposition is good.

\smallskip

Let $q\in k^{\times}$ and let $H$ be the algebra generated by $D_1, D_2, \sigma ^{\pm 1}$, subject to the relations
$$
D_1D_2 = D_2D_1,\quad \sigma\sigma^{-1} = \sigma ^{-1}\sigma =1\quad\text{and}\quad q\sigma D_i = D_i\sigma \quad\text{for $i=1,2$,}
$$
It is easy to check that $H$ is a Hopf algebra with
\begin{xalignat*}{3}
\Delta(D_1) &:= D_1\ot\sigma + 1\ot D_1, & \epsilon(D_1) &:= 0, & S(D_1) &:= -D_1\sigma ^{-1},\\
\Delta(D_2) &:= D_2\ot 1 + \sigma \ot D_2, &\epsilon(D_2) &:= 0, & S(D_2) &:= -\sigma ^{-1}D_2,\\
\Delta(\sigma )&:=\sigma\ot \sigma , & \epsilon(\sigma )&:=1, & S(\sigma )&:=\sigma^{-1}.
\end{xalignat*}
If $q$ is a primitive $l$-root of unity with $l\ge 2$, then the ideal $I$ of $H$ generated by $D_1^l$ and $D_2^l$ is a Hopf ideal. So, the quotient $H/I$ is also a Hopf algebra. Let
$$
H_q := \begin{cases}H/I &\text{if $q$ is a primitive $l$-root of unity with $l\ge 2$,}\\ H &\text{if $q=1$ or it is not a root of unity.}\end{cases}
$$
The Hopf algebra $H_q$ was considered in the paper~\cite{W}, where it was proved that
$$
\exp_q(tD_1\ot D_2)\! :=\! \begin{cases}\displaystyle{\sum_{i=0}^{l-1}}\frac{1}{(i)!_q} (tD_1\ot D_2)^i &\!\text{if $q$ is a primitive $l$-root of unity ($l\ge 2$),}\\
\displaystyle{\sum_{i=0}^{\infty}}\frac{1}{(i)_q!}(tD_1\!\ot\! D_2)^i &\!\text{if $q=1$ or it is not a root of unity,}
\end{cases}
$$
is an UDF based on $H_q$.

\smallskip

\subsection{Good transpositions of $\mathbf{H}_\mathbf{q}$ on an algebra}\label{subsect:
Some transp of H_q}
One of our main purposes in this paper is to construct formal deformation of algebras by using the UDF $\exp_q(tD_1\!\ot\! D_2)$. By Theorem~\ref{uso de UDF para construir deformaciones},  it will be sufficient to obtain examples of $H_q$-module algebras $(A,s)$, whose underlying transpositions $s$ satisfy
\begin{equation}\label{eq0}
(s\ot H_q)\xcirc (H_q\ot s)(D_1\ot D_2\ot a)= a\ot D_1\ot D_2\quad\text{for all $a\in A$.}
\end{equation}

\begin{definition} A $k$-linear map $s\colon H_q\ot A\to A\ot H_q$ is {\em good} if condition~\eqref{eq0} is fulfilled.
\end{definition}

It is evident that $s\colon H_q\ot A\to A\ot H_q$ is good if and only if there exists a bijective $k$-linear map $\alpha\colon A\to A$ such that
$$
s(D_1\ot a)= \alpha(a)\ot D_1\quad\text{and}\quad s(D_2\ot a)= \alpha^{-1}(a)\ot D_2\quad \text{for all $a\in A$.}
$$

\begin{lemma}\label{lema: las transposiciones de H_q sobre A inducen transposiciones de la subálgebra de grupo de H_q sobre A} Let $k[\sigma^{\pm 1}]$ denote the subHopfalgebra of $H_q$ generated by $\sigma$. Each transposition $s\colon H_q\ot A\to A\ot H_q$ takes $k[\sigma^{\pm 1}]\ot A$ onto $A\ot k[\sigma^{\pm 1}]$.
\end{lemma}

\begin{proof} Let $\tau$ be the flip. Since $\tau\xcirc s^{-1}\xcirc \tau$ is a transposition, it suffices to prove that $s(\sigma^{\pm 1}\ot a)\in A\ot k[\sigma^{\pm 1}]$ for all $a\in A$. Write
$$
s(\sigma\ot a)=\sum_{ijk}\gamma_{ijk}(a)\ot \sigma^iD_1^jD_2^k.
$$
Since $S^2(D_1) = q^{-1}D_1$, $S^2(D_2) = qD_2$ and $S^2(\sigma^{\pm 1}) = \sigma^{\pm 1}$, we have
\begin{align*}
\sum_{ijk}\gamma_{ijk}(a)\ot \sigma^iD_1^jD_2^k & = s(\sigma\ot a)\\[-7pt]
& = s\xcirc (S^2\ot A)(\sigma\ot a)\\
& = (A\ot S^2)\xcirc s(\sigma\ot a)\\
&= \sum_{ijk} q^{k-j} \gamma_{ijk}(a)\ot \sigma^iD_1^jD_2^k,
\end{align*}
and so $\gamma_{ijk}=0$ for $j\ne k$. Using now that
\begin{align*}
\sum_{ij}\gamma_{ijj}(a)\ot \Delta(\sigma)^i\Delta(D_1)^j\Delta(D_2)^j & = (A\ot \Delta)\xcirc s(\sigma\ot a)\\[-6pt]
& = (s\ot H_q)\xcirc (H_q\ot s)\xcirc (\Delta\ot A)(\sigma\ot a)\\
& = \sum_{iji'j'} \gamma_{i'j'j'}(\gamma_{ijj}(a))\ot \sigma^{i'}D_1^{j'}D_2^{j'}\ot \sigma^iD_1^jD_2^j,
\end{align*}
it is easy to check that $\gamma_{ijj}=0$ if $j>0$ (use that in each term of the right side the exponent of $D_1$ equals the exponent of $D_2$). For $\sigma^{-1}$ the same argument carries over. This finishes the proof.
\end{proof}

In the following result we obtain a characterization of the good transpositions of $H_q$ on an algebra $A$.

\begin{theorem}\label{algunas transposiciones de Hq} The following facts hold:

\begin{enumerate}

\smallskip

\item If $s\colon H_q\ot A\to A\ot H_q$ is a good transposition, then $s(\sigma^{\pm 1}\ot a)=a \ot \sigma^{\pm 1}$ for all $a\in A$ and the map $\alpha\colon A\to A$, defined by $s(D_1\ot a)= \alpha(a)\ot D_1$, is an algebra homomorphism.

\smallskip

\item Given an algebra automorphism $\alpha\colon A\to A$, there exists only one good transposition $s\colon H_q\ot A\to A\ot H_q$ such that $s(D_1\ot a)=\alpha (a)\ot D_1$ for all $a\in A$.

\smallskip

\end{enumerate}

\end{theorem}

\begin{proof} (1)\enspace  By Lemma~\ref{lema: las transposiciones de H_q sobre A inducen transposiciones de la subálgebra de grupo de H_q sobre A}, we know that $s$ induces by restriction a transposition of $k[\sigma^{\pm 1}]$ on $A$. Hence, by~\cite[Theorem~4.14]{G-G1}, there is a superalgebra structure $A= A_{+}\oplus A_{-}$ such that
$$
s(\sigma^i\ot a)=\begin{cases} a\ot \sigma^i & \text{if $a\in A_{+}$,}\\ a \ot \sigma^{-i} &\text{if $a\in A_{-}$.}\end{cases}
$$
Let $\alpha\colon A\to A$ be as in the statement. Since $\sigma$ is a transposition, if $a\in A_{-}$, then
\begin{align*}
\alpha(a)\ot D_1\ot \sigma + \alpha(a)\ot 1\ot D_1&=(A\ot \Delta) \xcirc s(D_1\ot a)\\
& = (s\ot H_q)\xcirc (H_q\ot s)\xcirc (\Delta \ot A )(D_1\ot a)\\
&= \alpha(a)\ot D_1\ot \sigma^{-1} + \alpha(a)\ot 1\ot D_1.
\end{align*}
So, $A_{-}=0$. Finally, $\alpha$ is an algebra homomorphism, because
$$
s(h\ot 1)\!=\!1\ot h\text{ for each $h\in H_q$}\!\quad\text{and}\!\quad s\xcirc (H_q\ot\mu_A)\!=\! (\mu_A\ot H_q)\xcirc (A\ot s)\xcirc (s\ot A).
$$

\smallskip

\noindent (2)\enspace By item~(1) and the comment preceding Lemma~\ref{lema: las transposiciones de H_q sobre A inducen transposiciones de la subálgebra de grupo de H_q sobre A}, it must be
$$
s(\sigma^{\pm 1}\ot a) = a\ot \sigma^{\pm 1},\quad\! s(D_1\ot a) = \alpha(a)\ot D_1\quad\! \text{and}\quad \!s(D_2\ot a) = \alpha^{-1}(a)\ot D_2.
$$
So, necessarily
$$
s(\sigma^iD_1^jD_2^k\ot a) = \alpha^{j-k}(a)\ot \sigma^iD_1^jD_2^k.
$$
We leave to the reader the task to prove that $s$ is a good transposition.
\end{proof}

\subsection{Some $\mathbf{H}_\mathbf{q}$-module algebra structures}\label{subsection: some Hq-module algebra structures}
Let $A$ be an algebra. Let us consider $k$-linear maps $\varsigma, \delta_1, \delta_2 \colon A\to A$. It is evident that there is a (necessarily unique) action $\rho\colon H_q\ot A\to A$ such that
\begin{equation}
\rho(\sigma \ot a) = \varsigma(a),\quad\rho(D_1\ot a) = \delta_1(a)\quad\text{and} \quad\rho(D_2\ot a) = \delta_2(a)\label{eq1}
\end{equation}
for all $a\in A$, if and only if the maps $\varsigma$, $\delta_1$ and $\delta_2$ satisfy the following conditions:
\begin{enumerate}

\smallskip

\item $\varsigma$ is a bijective map,

\smallskip

\item $\delta_1\xcirc\delta_2 = \delta_2\xcirc\delta_1$,

\smallskip

\item $q\varsigma\xcirc\delta_i = \delta_i\xcirc\varsigma$ for $i=1,2$.

\smallskip

\item If $q\ne 1$ and $q^l = 1$, then $\delta_1^l = \delta_2^l=0$.

\smallskip

\end{enumerate}

Let $s\colon H_q\ot A\to A\ot H_q$ be a good transposition and let $\alpha$ be the associated automorphism. Let $\varsigma$, $\delta_1$ and $\delta_2$ be $k$-linear endomorphisms of $A$ satisfying (1)--(4). Next, we determine the conditions that $\varsigma$, $\delta_1$ and $\delta_2$ must satisfy in order that $(A,s)$ becomes an $H_q$-module algebra via the action $\rho$ defined by~\eqref{eq1}.

\begin{theorem}\label{algunas acciones de H_q sobre algebras} $(A,s)$ is an $H_q$-module algebra via $\rho$ if and only if:

\begin{enumerate}

\smallskip

\item[(5)] $\varsigma$ is an algebra automorphism,

\smallskip

\item[(6)] $\alpha\xcirc\delta_i = \delta_i\xcirc\alpha$ for $i=1,2$,

\smallskip

\item[(7)] $\alpha\xcirc\varsigma = \varsigma\xcirc\alpha$,

\smallskip

\item[(8)] $\delta_i(1)=0$ for $i=1,2$,

\smallskip

\item[(9)] $\delta_1(ab) = \delta_1(a)\varsigma(b) + \alpha(a)\delta_1(b)$ for all $a,b\in A$,

\smallskip

\item[(10)] $\delta_2(ab) = \delta_2(a)b + \varsigma\bigl(\alpha^{-1}(a)\bigr)\delta_2(b)$ for all $a,b\in A$.

\smallskip

\end{enumerate}
\end{theorem}

\begin{proof} Assume that $(A,s)$ is an $H_q$-module algebra and let $\tau\colon H_q\ot H_q\to H_q\ot H_q$ be the flip. Evaluating the equality
$$
s\xcirc (H_q\ot\rho) = (\rho\ot H_q)\xcirc (H_q\ot s)\xcirc (\tau\ot A)
$$
successively on $D_1\ot D_i\ot a$ and $D_1\ot \sigma\ot a$ with $i\in \{1,2\}$ and $a\in A$ arbitrary, we verify that items~(6) and~(7) are satisfied. Item~(8) follows from the fact that $D_1\!\cdot\!1 = D_2\!\cdot\!1 = 0$. Finally, using that $\sigma\!\cdot\!1 = 1$ and evaluating the equality
$$
\rho\xcirc(H_q\ot\mu_A) = \mu_A\xcirc(\rho\ot\rho)\xcirc(H_q\ot s\ot A)\xcirc (\Delta\ot A\ot A)
$$
on $\sigma \ot a\ot b$ and $D_i\ot a\ot b$, with $i=1,2$ and $a,b\in A$ arbitrary, we see that items~(5), (9) and~(10) hold. So, conditions~(5)--(10) are necessary. By Remark~\ref{condiciones para que (A,s) sea un H-módulo álgebra a izquierda}, in order to verify that they are also sufficient, it is enough to check that they imply that
\begin{align*}
& h\!\cdot\!1 = \epsilon (h),\\
&s(h\ot l\!\cdot\!a) = (\rho\ot H_q)\xcirc (H_q\ot s)(l\ot h\ot a),\\
& h\cdot(ab) = \mu_A\xcirc(\rho\ot\rho)\xcirc(H_q\ot s\ot A)(\Delta(h)\ot a\ot b),
\end{align*}
for all $a,b\in A$ and $h,l\in \{D_1,D_2,\sigma^{\pm 1}\}$. We leave this task to the reader.
\end{proof}

Note that condition~(8) in Theorem~\ref{algunas acciones de H_q sobre algebras} is redundant since it can be obtained by applying condition~(9) and~(10) with $a=b=1$.

\section{$\mathbf{H}_\mathbf{q}$-module algebra structures on crossed products}\label{Hq-module algebra structures on crossed products}
Let $G$ be a group endowed with a representation on a $k$-vector space $V$ of dimension~$n$. Consider the symmetric $k$-algebra $S(V)$ equipped with the unique action of $G$ by automorphisms that extends the action of $G$ on $V$  and take $A = S(V)\#_f G$, where $f\colon G\times G\to k^{\times}$ is a normal cocycle. By definition the $k$-algebra $A$ is a free left $S(V)$-module with basis $\{w_g:g\in G\}$. Its product is given by
$$
(Pw_g)(Qw_h) :=P{}^g\!Qf(g,h)w_{gh},
$$
where ${}^g\! Q$ denotes the action of $g$ on $Q$. This section is devoted to the study of the $H_q$-module algebras $(A,s)$, with $s$ good, that satisfy:
$$
s(H_q\ot V)\subseteq V\ot H_q,\!\quad s(H_q\ot kw_g)\subseteq kw_g\ot H_q,\!\quad \sigma\!\cdot\! v\in V\!\quad\text{and}\!\quad \sigma\!\cdot\!w_g\in kw_g,
$$
for all $v\in V$ and $g\in G$. In Theorem~\ref{condiciones para la existencia de una estructura de H_q modulo algebra sobre (A,s) cuando A es un producto cruzado} we give a general characterization of these module algebras, and in Subsection~\ref{First case} we consider a specific case which is more suitable for finding concrete examples, and we study it in detail. Finally in Subsection~\ref{Second case} we consider the case where the cocycle involves several non necessarily central elements of $G$.

\smallskip

In the following proposition we characterize the good transpositions $s$ of $H_q$ on $A$  satisfying the hypothesis mentioned above. By theorem~\ref{algunas transposiciones de Hq} this is equivalent to require that the $k$-linear map $\alpha\colon A\to A$ associated with $\alpha$, takes $V$ to $V$ and $kw_g$ to $kw_g$ for all $g\in G$.

\begin{proposition}\label{prop: algunas transposiciones buenas} Let $\hat{\alpha}\colon V\to V$ be a $k$-linear map and $\chi_{\alpha}\colon G\to k^{\times}$ a map. There is a good transposition $s\colon H_q\ot A\to A\ot H_q$, such that
$$
s(D_1\ot v)=\hat{\alpha}(v)\ot D_1\quad\text{and}\quad s(D_1\ot w_g) = \chi_{\alpha}(g)w_g\ot D_1
$$
for all $v\in V$ and $g\in G$, if and only if $\hat{\alpha}$ is a bijective $k[G]$-linear map and $\chi_{\alpha}$ is a group homomorphism.
\end{proposition}

\begin{proof} By Theorem~\ref{algunas transposiciones de Hq} we know that $s$ exists if an only  if the $k$-linear map \hbox{$\alpha\colon\! A\to A$} defined by
$$
\alpha(v_1\cdots v_mw_g) := \hat{\alpha}(v_1)\cdots\hat{\alpha}(v_m)\chi_{\alpha}(g)w_g,
$$
is an automorphism. But, if this happens, then

\smallskip

\noindent a)\enspace $\chi_{\alpha}$ is a morphism since
$$
\chi_{\alpha}(g)\chi_{\alpha}(h)f(g,h)w_{gh}=\alpha(w_g)\alpha(w_h)=\alpha(w_gw_h)= \chi_{\alpha}(gh)f(g,h)w_{gh}
$$
for all $g,h\in G$,

\smallskip

\noindent b)\enspace $\hat{\alpha}$ is a bijective $k[G]$-linear map, since it is the restriction and correstriction of $\alpha$ to $V$, and
$$
\hat{\alpha}({}^g\!v) = \alpha(w_g)\hat{\alpha}(v)\alpha(w_g^{-1}) = \chi_{\alpha}(g)w_g \hat{\alpha}(v)(\chi_{\alpha}(g)w_g)^{-1} = w_g\hat{\alpha}(v)w_g^{-1} = {}^g\!\hat{\alpha}(v).
$$
Conversely, if $\hat{\alpha}$ is a bijective map then $\alpha$ is also, and if $\hat{\alpha}$ is a $k[G]$-linear map and $\chi_{\alpha}$ is a morphism, then
$$
\alpha(w_g)\hat{\alpha}(v) = \chi_{\alpha}(g)w_g\hat{\alpha}(v)\! = {}^g\! \hat{\alpha}(v)\chi_{\alpha}(g)w_g\! =\hat{\alpha}({}^g\!v)\alpha(w_g)
$$
and
$$
\alpha(w_g)\alpha(w_h) = \chi_{\alpha}(g)w_g \chi_{\alpha}(h)w_h\! = f(g,h)\chi_{\alpha}(gh) w_{gh}\! =\alpha(f(g,h)w_{gh}),
$$
for all $v\in V$ and $g,h\in G$, from which it follows easily that $\alpha$ is a morphism.
\end{proof}

Let $A = S(V)\#_f G$ be as above. Throughout this section we fix a morphism $\chi_{\alpha} \colon G\to k^{\times}$ and a bijective $k[G]$-linear map $\hat{\alpha}\colon V\to V$, and we let $\alpha\colon A\to A$ denote the automorphism determined by $\hat{\alpha}$ and $\chi_{\alpha}$. Moreover we will call
$$
s\colon H_q\ot A \to A\ot H_q
$$
the good transposition associated with $\alpha$. Our purpose is to obtain all the $H_q$-mo\-dule algebra structures on $(A,s)$ such that
\begin{equation}
\sigma\!\cdot\!v\in V\quad\text{and}\quad \sigma\!\cdot\!w_g\in kw_g\quad\text{for all $v\in V$ and $g\in G$.}\label{eq2}
\end{equation}
Under these restrictions we obtain conditions which allow us to construct all $H_q$-mo\-dule structures in concrete examples. Thanks to Theorem~\ref{uso de UDF para construir deformaciones} and the fact that $\exp_q(tD_1\ot D_2)$ is an UDF based on $H_q$, each one of these examples produces automatically a formal deformation of $A$. First note that given an $H_q$-module algebra structure on $(A,s)$ satisfying~\eqref{eq2}, we can define $k$-linear maps
$$
\hat{\delta}_1\colon V\to A,\quad\hat{\delta}_2\colon V\to A\quad\text{and}\quad \hat{\varsigma}\colon V\to V
$$
and maps
$$
\ov{\delta}_1\colon G\to A,\quad \ov{\delta}_2\colon G\to A\quad\text{and}\quad \chi_{\varsigma} \colon G\to k^{\times},
$$
by
$$
\hat{\delta}_i(v) := D_i\!\cdot\!v,\quad \hat{\varsigma}(v) := \sigma\!\cdot\!v,\quad \ov{\delta}_i(g) := D_i\!\cdot\!w_g \quad \text{and}\quad \sigma\!\cdot\!w_g := \chi_{\varsigma}(g)w_g.
$$

\begin{lemma}\label{condicion para que varsigma sea un automorfismo de algebras} Let $\hat{\varsigma}\colon V\to V$ be a $k$-linear map and $\chi_{\varsigma}\colon G\to k^{\times}$ be a map. Then, the map $\varsigma\colon A \to A$ defined by
$$
\varsigma(\mathbf{v}_{1m}w_g) := \hat{\varsigma}(v_1)\cdots\hat{\varsigma}(v_m) \chi_{\varsigma}(g) w_g,
$$
is a $k$-algebra automorphism if and only if $\hat{\varsigma}$ is a bijective $k[G]$-linear map and $\chi_{\varsigma}$ is a group homomorphism.
\end{lemma}

\begin{proof} This was checked in the proof of Proposition~\ref{prop: algunas transposiciones buenas}.
\end{proof}

\begin{lemma}\label{buena def de delta_i} Let $\hat{\delta}_1\colon V\to A$ and $\hat{\delta}_2\colon V\to A$ be $k$-linear maps and let $\ov{\delta}_1\colon G\to A$ and $\ov{\delta}_2\colon G\to A$ be maps.

\begin{enumerate}

\smallskip

\item The $k$-linear map $\delta_1\colon A\to A$ given by
$$
\qquad  \delta_1(\mathbf{v}_{1m} w_g) := \sum_{j=1}^m \alpha(\mathbf{v}_{1,j-1}) \hat{\delta}_1(v_j) \varsigma(\mathbf{v}_{j+1,m}w_g) + \alpha(\mathbf{v}_{1m}) \ov{\delta}_1(g),
$$
where $\mathbf{v}_{hl} = v_h\cdots v_l$, is well defined if and only if
\begin{equation}
\qquad\quad\,\,\,\, \hat{\delta}_1(v)\hat{\varsigma}(w) + \hat{\alpha}(v) \hat{\delta}_1(w) = \hat{\delta}_1(w) \hat{\varsigma}(v) + \hat{\alpha}(w) \hat{\delta}_1(v)\quad\text{for all $v,w\in V$.}\label{eq2.1}
\end{equation}

\smallskip

\item The map $\delta_2\colon A\to A$ given by
$$
\qquad \delta_2(\mathbf{v}_{1m} w_g) := \sum_{j=1}^m\varsigma\bigl(\alpha^{-1} (\mathbf{v}_{1,j-1})\bigr) \hat{\delta}_2(v_j)\mathbf{v}_{j+1,m}w_g + \varsigma\bigl(\alpha^{-1}\bigr)(\mathbf{v}_{1m}) \ov{\delta}_2(g)
$$
is well defined if and only if
\begin{equation}
\qquad\quad\,\,\,\,\hat{\delta}_2(v) w\!+\! \varsigma\bigl(\!{\hat{\alpha}}^{-1}(v)\!\bigr) \hat{\delta}_2(w) = \hat{\delta}_2(w) v\! +\! \varsigma\bigl(\!{\hat{\alpha}}^{-1}(w)\! \bigr) \hat{\delta}_2(v)\!\!\quad\text{for all $v,w\in V$.}\label{eq2.2}
\end{equation}

\end{enumerate}

\end{lemma}

\begin{proof} We prove the first assertion and leave the second one, which is similar, to the reader. The only if part follows immediately by noting that
$$
\hat{\delta}_1(v)\hat{\varsigma}(w) + \hat{\alpha}(v) \hat{\delta}_1(w)=\delta_1(vw)=\delta_1(wv) = \hat{\delta}_1(w) \hat{\varsigma}(v) + \hat{\alpha}(w) \hat{\delta}_1(v).
$$
In order to prove the if part it suffices to check that
$$
\delta_1(v_1\cdots v_{i-1}v_{i+1}v_iv_{i+2}\cdots v_mw_g)=\delta_1(\mathbf{v}_{1m}w_g)\quad \text{for all $i<m$,}
$$
which follows easily from the hypothesis.
\end{proof}

\begin{lemma}\label{propiedades de delta_i} Assume that $\varsigma$ is an algebra automorphism and $\delta_1$, $\delta_2$ are well defined. The following facts hold:

\begin{enumerate}

\smallskip

\item The map $\delta_1$ satisfies
$$
\delta_1(x_1\cdots x_m) = \sum_{j=1}^m\alpha(x_1\dots x_{j-1})\delta_1(x_j) \varsigma(x_{j+1} \cdots x_m)
$$
for all $x_1,\dots, x_m\in k\#_fG \cup V$, if and only if

\begin{enumerate}

\smallskip

\item $\hat{\delta}_1({}^g\! v)\chi_{\varsigma}(g)w_g + \hat{\alpha}({}^g\! v) \ov{\delta}_1(g) = \ov{\delta}_1(g) \hat{\varsigma}(v) + \chi_{\alpha}(g)w_g \hat{\delta}_1(v)$,

\smallskip

\item $f(g,h)\ov{\delta}_1(gh) = \ov{\delta}_1(g)\chi_{\varsigma}(h)w_h + \chi_{\alpha}(g)w_g\ov{\delta}_1(h)$,

\smallskip

\end{enumerate}
for all $v\in V$ and $g,h\in G$.

\smallskip

\item The map $\delta_2$ satisfies
$$
\delta_2(x_1\cdots x_m) = \sum_{j=1}^m\varsigma\xcirc \alpha^{-1}(x_1\dots x_{j-1}) \delta_1(x_j) x_{j+1}\cdots x_m
$$
for all $x_1,\dots,x_m\in k\#_fG \cup V$, if and only if

\begin{enumerate}

\smallskip

\item $\hat{\delta}_2({}^g\! v)w_g + \hat{\varsigma}\bigl({\hat{\alpha}}^{-1} ({}^g\! v)\bigr) \ov{\delta}_2(g) = \ov{\delta}_2(g) v + \chi_{\varsigma}(g) \chi_{\alpha}^{-1}(g) w_g \hat{\delta}_2(v)$,

\smallskip

\item $f(g,h)\ov{\delta}_2(gh) = \ov{\delta}_2(g)w_h + \chi_{\varsigma}(g) \chi_{\alpha}^{-1}(g) w_g \ov{\delta}_2(h)$,

\smallskip

\end{enumerate}
for all $v\in V$ and $g,h\in G$.
\end{enumerate}
\end{lemma}

\begin{proof} We prove the first assertion and leave the second one to the reader. For the only if part it suffices to note that
\begin{align*}
& \hat{\delta}_1({}^g\! v)\varsigma(w_g) + \alpha({}^g\! v)\ov{\delta}_1(g) = \delta_1({}^g\! vw_g) = \delta_1(w_g v) = \ov{\delta}_1(g) \varsigma(v) + \alpha(w_g) \hat{\delta}_1(v),\\
& f(g,h)\ov{\delta}_1(gh) =\delta_1(w_gw_h)= \ov{\delta}_1(g)\varsigma(w_h) + \alpha(w_g) \ov{\delta}_1(h),
\end{align*}
and to use the definitions of $\varsigma(w_g)$ and $\alpha(w_g)$. We prove the sufficient part by induction on $r=m+1-i$, where $i$ is the first index with $x_i\in k\#_f G$ (if \hbox{$x_1,\dots, x_m\in V$} we set $r := 0$). For $r\in \{0,1\}$ the result follows immediately from the definition of $\delta_1$. Assume that it is true when $r<r_0$ and that $m+1-i = r_0$. If $x_i = w_g$ and $x_{i+1} = v\in V$, then
$$
\delta_1(x_1\cdots x_m) = \delta_1(y_1\cdots y_m)\quad\text{where }  y_j = \begin{cases} x_j & \text{if $j\notin \{i,i+1\}$,}\\ {}^g\! v & \text{if $j = i$,} \\ w_g & \text{if $j = i+1$,}\end{cases}
$$
and hence, by the inductive hypothesis and item~(a),
\begin{align*}
\delta_1(x_1\cdots x_m) & = \sum_{j=1}^m\alpha(y_1\dots y_{j-1})\delta_1(y_j) \varsigma(y_{j+1} \cdots y_m)\\
& = \sum_{j=1}^m \alpha(x_1\dots x_{j-1})\delta_1(x_j)\varsigma(x_{j+1}\cdots x_m).
\end{align*}
If $x_i = w_g$ and $x_{i+1} = w_h$, then
$$
\delta_1(x_1\cdots x_m) = f(g,h)\delta_1(y_1\cdots y_{m-1})\quad\text{where }  y_j = \begin{cases} x_j & \text{if $j<i$,}\\ w_{gh} & \text{if $j = i$,}\\ x_{j+1}& \text{if $j>i$,}\end{cases}
$$
and hence, by the inductive hypothesis and item~(b),
\begin{align*}
\delta_1(x_1\cdots x_m) & = \sum_{j=1}^{m-1} f(g,h) \alpha(y_1\dots y_{j-1})\delta_1(y_j) \varsigma(y_{j+1}\cdots y_{m-1})\\
& = \sum_{j=1}^m \alpha(x_1\dots x_{j-1})\delta_1(x_j) \varsigma(x_{j+1}\cdots x_m),
\end{align*}
as we want.
\end{proof}

\begin{theorem}\label{condiciones para la existencia de una estructura de H_q modulo algebra sobre (A,s) cuando A es un producto cruzado} Let $\hat{\delta}_1\colon V\to A$, $\hat{\delta}_2\colon V\to A$ and $\hat{\varsigma}\colon V\to V$ be $k$-linear maps and let $\ov{\delta}_1\colon G\to A$, $\ov{\delta}_2\colon G\to A$ and  $\chi_{\varsigma}\colon G\to k^{\times}$ be maps. There is an $H_q$-module algebra structure on $(A,s)$, such that
$$
\sigma \!\cdot\!v = \hat{\varsigma}(v),\quad \sigma \!\cdot\!w_g=\chi_{\varsigma}(g)w_g, \quad D_i \!\cdot\!v = \hat{\delta}_i(v) \quad\text{and}\quad D_i \!\cdot\!w_g = \ov{\delta}_i(g)
$$
for all $v\in V$, $g\in G$ and $i\in \{1,2\}$, if and only if

\begin{enumerate}

\smallskip

\item $\hat{\varsigma}\colon V\to V$ is a bijective $k[G]$-linear map and $\chi_{\varsigma}$ is a group homomorphism,

\smallskip

\item Conditions~\eqref{eq2.1} and~\eqref{eq2.2} in Lemma~\ref{buena def de delta_i} and items~(1)(a), (1)(b), (2)(a) and~(2)(b) in Lemma~\ref{propiedades de delta_i} are satisfied,

\smallskip

\item $\hat{\delta}_i\xcirc \hat{\alpha} = \alpha\xcirc \hat{\delta}_i$,

\smallskip

\item $\chi_{\alpha}(g) \ov{\delta}_i(g) = \alpha\bigl(\ov{\delta}_i(g)\bigr)$ for all $g\in G$,

\smallskip

\item $\hat{\varsigma}\xcirc \hat{\alpha} = \hat{\alpha}\xcirc\hat{\varsigma}$

\smallskip

\item The maps $\varsigma\colon A\to A$, $\delta_1\colon A\to A$ and $\delta_2\colon A\to A$, introduced in Lemmas~\ref{condicion para que varsigma sea un automorfismo de algebras} and~\ref{buena def de delta_i}, satisfy the following properties:
\begin{xalignat*}{3}
&\qquad\quad\delta_2\xcirc \hat{\delta}_1= \delta_1\xcirc \hat{\delta}_2,&& \hat{\delta}_i\xcirc \hat{\varsigma} = q \varsigma \xcirc\hat{\delta}_i,&& \delta_2\xcirc \ov{\delta}_1 = \delta_1\xcirc\ov{\delta}_2,\\
&\qquad\quad \chi_{\varsigma}(g) \ov{\delta}_i(g) = q \varsigma\bigl(\ov{\delta}_i(g) \bigr), && \delta_1^l=\delta_2^l=0\text{ if $q\ne 1$ and $q^l=1$.}
\end{xalignat*}

\end{enumerate}
\end{theorem}

\begin{proof} By Theorem~\ref{algunas acciones de H_q sobre algebras} and the discussion above it, we know that to have an \hbox{$H_q$-module} algebra structure on $(A,s)$ satisfying the requirements in the statement is equivalent to have maps \mbox{$\varsigma,\delta_1,\delta_2 \colon A\to A$} satisfying conditions~(1)--(10) in Subsection~\ref{subsection: some Hq-module algebra structures} and such that
$$
\varsigma(v) = \hat{\varsigma}(v),\quad \varsigma(w_g)=\chi_{\varsigma}(g)w_g, \quad \delta_i(v) = \hat{\delta}_i(v) \quad\text{and}\quad \delta_i(w_g) = \ov{\delta}_i(g)
$$
for all $v\in V$, $g\in G$ and $i\in \{1,2\}$. Now, it is easy to see that

\smallskip

\noindent a)\enspace If $\varsigma$, $\delta_1$ and $\delta_2$ satisfy conditions~(5), (9) and (10) in Subsection~\ref{subsection: some Hq-module algebra structures}, then
\begin{align*}
&\varsigma(\mathbf{v}_{1m}w_g) = \hat{\varsigma}(v_1)\cdots\hat{\varsigma}(v_m) \chi_{\varsigma}(g) w_g,\\
&\delta_1(\mathbf{v}_{1m} w_g) = \sum_{j=1}^m \alpha(\mathbf{v}_{1,j-1}) \hat{\delta}_1(v_j) \varsigma(\mathbf{v}_{j+1,m}w_g) + \alpha(\mathbf{v}_{1m}) \ov{\delta}_1(g),\\
& \delta_2(\mathbf{v}_{1m} w_g) = \sum_{j=1}^m \varsigma\bigl(\alpha^{-1}(\mathbf{v}_{1,j-1}) \bigr) \hat{\delta}_2(v_j) \mathbf{v}_{j+1,m}w_g + \varsigma\bigl(\alpha^{-1}(\mathbf{v}_{1m}) \bigr) \ov{\delta}_2(g),
\end{align*}
where $\mathbf{v}_{hl} = v_h\cdots v_l$.

\smallskip

\noindent b)\enspace By Lemmas~\ref{condicion para que varsigma sea un automorfismo de algebras}, \ref{buena def de delta_i} and~\ref{propiedades de delta_i}, the maps defined in~a) satisfy conditions~(1), (5), (8), (9) and~(10) in Subsection~\ref{subsection: some Hq-module algebra structures} if and only if items~(1) and~(2) of the present theorem are fulfilled.

\smallskip

\noindent So, in order to finish the proof it suffices to check that:

\smallskip

\noindent c)\enspace Conditions~(6) and~(7) in Subsection~\ref{subsection: some Hq-module algebra structures} are satisfied if and only if items~(3)--(5) of the present theorem are fulfilled,

\smallskip

\noindent d)\enspace Conditions~(2), (3) and~(4) in Subsection~\ref{subsection: some Hq-module algebra structures} are satisfied if and only if item~(6) of the present theorem is fulfilled.

\smallskip

\noindent We leave this task to the reader.
\end{proof}

We are going now to consider several particular cases, with the purpose of obtaining more precise results. This will allow us to give some specific examples of formal deformations of associative algebras.

\subsection{First case}\label{First case} Let $\hat{\alpha}$, $\chi_{\alpha}$, $\alpha$ and $s$ be as in the discussion following Proposition~\ref{prop: algunas transposiciones buenas}. Let $\hat{\delta}_1\colon V\to A$, $\hat{\delta}_2\colon V\to A$ and $\hat{\varsigma}\colon V\to V$ be $k$-linear maps and let  $\chi_{\varsigma}\colon G\to k^{\times}$ be a map. Assume that the kernels of $\hat{\delta}_1$ and $\hat{\delta}_2$ have codimension~$1$, $\ker\hat{\delta}_1\ne \ker\hat{\delta}_2$ and there exist $x_i\in V\setminus\ker\hat{\delta}_i$, such that $\hat{\delta}_i(x_i) = P_iw_{g_i}$ with $P_i\in S(V)$ and $g_i\in G$. Without loss of generality we can assume that $x_1\in \ker\hat{\delta}_2$ and $x_2\in \ker\hat{\delta}_1$ (and we do it). For $g\in G$ and $i\in\{1,2\}$, let $\lambda_{ig},\omega_i,\nu_i\in k$ be the elements defined by the following conditions:
$$
{}^g\! x_i - \lambda_{ig} x_i\in \ker \hat{\delta}_i,\quad \hat{\varsigma}(x_i) - \omega_i x_i\in \ker \hat{\delta}_i\quad\text{and}\quad \hat{\alpha}(x_i) - \nu_i x_i\in \ker\hat{\delta}_i.
$$

\begin{theorem}\label{condiciones para la existencia de estructuras de H_q-modulo algebra (caso 1)} There is an $H_q$-module algebra structure on $(A,s)$, satisfying
$$
\sigma \!\cdot\!v = \hat{\varsigma}(v),\quad \sigma \!\cdot\!w_g=\chi_{\varsigma}(g)w_g, \quad D_i \!\cdot\!v = \hat{\delta}_i(v) \quad\text{and}\quad D_i \!\cdot\!w_g =0
$$
for all $v\in V$, $g\in G$ and $i\in \{1,2\}$,  if and only if

\begin{enumerate}

\smallskip

\item $\hat{\varsigma}$ is a bijective $k[G]$-linear map and $\chi_{\varsigma}$ is a group homomorphism,

\smallskip

\item $\hat{\varsigma}(v) = {}^{g_1^{-1}}\!\hat{\alpha}(v)$ for all $v\in \ker \hat{\delta}_1$ and $\hat{\varsigma}(v) = {}^{g_2}\!\hat{\alpha}(v)$ for all $v\in \ker \hat{\delta}_2$,

\smallskip

\item $g_1$ and $g_2$ belong to the center of $G$,

\smallskip

\item $\ker \hat{\delta}_1$ and $\ker \hat{\delta}_2$ are $G$-submodules of $V$,

\smallskip

\item ${}^g\! P_1 = \lambda_{1g} \chi_{\alpha}^{-1}(g)\chi_{\varsigma}(g)f^{-1}(g,g_1) f(g_1,g) P_1$ for all $g\in G$,

\smallskip

\item ${}^g\! P_2 = \lambda_{2g}\chi_{\alpha}(g)\chi_{\varsigma}^{-1}(g)f^{-1}(g,g_2) f(g_2,g)P_2$ for all $g\in G$,

\smallskip

\item $\hat{\alpha}(\ker \hat{\delta}_i) = \ker \hat{\delta}_i$ for $i\in \{1,2\}$,

\smallskip

\item $P_1\in \ker\delta_2$ and $P_2\in \ker\delta_1$, where $\delta_1$ and $\delta_2$ are the maps defined by
\begin{align*}
&\qquad\delta_1(\mathbf{v}_{1m} w_g) := \sum_{j=1}^m \alpha(\mathbf{v}_{1,j-1}) \hat{\delta}_1(v_j) \varsigma(\mathbf{v}_{j+1,m}w_g),\\
&\qquad\delta_2(\mathbf{v}_{1m} w_g) := \sum_{j=1}^m\varsigma\bigl(\alpha^{-1} (\mathbf{v}_{1,j-1})\bigr) \hat{\delta}_2(v_j)\mathbf{v}_{j+1,m}w_g,
\end{align*}
in which $\mathbf{v}_{hl} = v_h\cdots v_l$,

\smallskip

\item $\varsigma(P_i) = q^{-1} \omega_i\chi_{\varsigma}^{-1}(g_i)P_i$ and $\alpha(P_i) = \nu_i\chi_{\alpha}^{-1}(g_i)P_i$ for $i\in \{1,2\}$, where  $\varsigma$ is the map
    given by
$$
\varsigma(\mathbf{v}_{1m}w_g) = \hat{\varsigma}(v_1)\cdots\hat{\varsigma}(v_m) \chi_{\varsigma}(g) w_g,
$$

\smallskip

\item If $q\ne 1$ and $q^l=1$, then $\delta_1^l=\delta_2^l=0$.

\smallskip

\end{enumerate}

\end{theorem}

In order to prove this result we first need to establish some auxiliary results.

\begin{lemma}\label{buena def de delta_i caso 1} The following facts hold:

\begin{enumerate}

\smallskip

\item Condition~\eqref{eq2.1} of Lemma~\ref{buena def de delta_i} is satisfied if and only if ${}^{g_1}\!\hat{\varsigma}(v) = \hat{\alpha}(v)$ for all $v\in \ker \hat{\delta}_1$.

\smallskip

\item Condition~\eqref{eq2.2} of Lemma~\ref{buena def de delta_i} is satisfied if and only if ${}^{g_2}\! v = \hat{\varsigma}\bigl(\hat{\alpha}^{-1}(v)\bigr)$ for all $v\in \ker \hat{\delta}_2$.

\end{enumerate}

\end{lemma}

\begin{proof} We prove item~1) and we leave item~2), which is similar, to the reader. We must check that
\begin{equation}
\hat{\delta}_1(v)\hat{\varsigma}(w) + \hat{\alpha}(v) \hat{\delta}_1(w) = \hat{\delta}_1(w)\hat{\varsigma}(v) + \hat{\alpha}(w) \hat{\delta}_1(v)\quad\text{for all $v,w\in V$}\label{eqq1}
\end{equation}
if and only if $\hat{\varsigma}_1(v) = {}^{g_1^{-1}}\!{\hat{\alpha}}(v)$ for all $v\in \ker \hat{\delta}_1$. It is clear that we can suppose that $v,w\in \{x_1\}\cup \ker\hat{\delta}_1$.  When $v,w\in \ker\hat{\delta}_1$ or $v = w = x_1$ the equality~\eqref{eqq1} is trivial. Assume $v = x_1$ and $w\in \ker\hat{\delta}_1$. Then,
$$
\hat{\delta}_1(v)\hat{\varsigma}(w) + \hat{\alpha}(v) \hat{\delta}_1(w) = P_1 w_{g_1} \hat{\varsigma}(w) = P_1  {}^{g_1}\! \hat{\varsigma}(w) w_{g_1}
$$
and
$$
\hat{\delta}_1(w)\hat{\varsigma}(v) + \hat{\alpha}(w) \hat{\delta}_1(v) = \hat{\alpha}(w)P_1 w_{g_1} = P_1 \hat{\alpha}(w) w_{g_1}
$$
So, in this case, the result is true. Case $v \in \ker\hat{\delta}_1$ and $w = x_1$ can be treated in a similar way.
\end{proof}

\begin{lemma}\label{propiedades de delta_i caso 1} The following facts hold:

\begin{enumerate}

\smallskip

\item Items~(1)(a) and~(1)(b) of Lemma~\ref{propiedades de delta_i} are satisfied if and only if

\begin{enumerate}

\smallskip

\item $\ker \hat{\delta}_1$ is a $G$-submodule of $V$,

\smallskip

\item $g_1$ belongs to the center of $G$,

\smallskip

\item ${}^g\! P_1 = \lambda_{1g} \chi_{\alpha}^{-1}(g)\chi_{\varsigma}(g)f^{-1}(g,g_1) f(g_1,g) P_1$, for all $g\in G$.

\smallskip

\end{enumerate}

\smallskip

\item Items~(2)(a) and~(2)(b) of Lemma~\ref{propiedades de delta_i} are satisfied if and only if

\begin{enumerate}

\smallskip

\item $\ker \hat{\delta}_2$ is a $G$-submodule of $V$,

\smallskip

\item $g_2$ belongs to the center of $G$,

\smallskip

\item ${}^g\! P_2 = \lambda_{2g}\chi_{\alpha}(g)\chi_{\varsigma}^{-1}(g)f^{-1}(g,g_2) f(g_2,g)P_2$, for all $g\in G$.

\smallskip

\end{enumerate}
\end{enumerate}
\end{lemma}

\begin{proof} We prove item~1) and we leave item~2) to the reader. Since $\ov{\delta}_1 = 0$, it is sufficient to prove that
\begin{equation}
\hat{\delta}_1({}^g\! v)\chi_{\varsigma}(g)w_g = \chi_{\alpha}(g)w_g \hat{\delta}_1(v) \quad \text{for all $v\in V$ and $g\in G$,}\label{eqq2}
\end{equation}
if and only if conditions~(1)(a), (1)(b) and~(1)(c) are satisfied. We can assume that $v\in \{x_1\}\cup \ker\hat{\delta}_1$. When $v\in \ker\hat{\delta}_1$, then equality~\eqref{eqq2} is true if and only if ${}^g\! v\in \ker\hat{\delta}_1$. Now, since
$$
\hat{\delta}_1({}^g\! x_1)\chi_{\varsigma}(g)w_g = \lambda_{1g} P_1w_{g_1}\chi_{\varsigma}(g) w_g = \lambda_{1g} P_1\chi_{\varsigma}(g)f(g_1,g)w_{g_1g}
$$
and
$$
\chi_{\alpha}(g)w_g \hat{\delta}_1(x_1) = \chi_{\alpha}(g) w_g P_1w_{g_1} = \chi_{\alpha}(g) {}^g\! P_1 f(g,g_1)w_{gg_1},
$$
equality~\eqref{eqq2} is true for $v = x_1$ and $g\in G$ if and only if conditions~(1)(b) and~ (1)(c) are satisfied.
\end{proof}

\begin{proof}[Proof of Theorem~\ref{condiciones para la existencia de estructuras de H_q-modulo algebra (caso 1)}] First note that item~(1) coincide with item~(1) of Theorem~\ref{condiciones para la existencia de una estructura de H_q modulo algebra sobre (A,s) cuando A es un producto cruzado} and that,  by Lemmas~\ref{buena def de delta_i caso 1} and~\ref{propiedades de delta_i caso 1}, item~(2) of Theorem~\ref{condiciones para la existencia de una estructura de H_q modulo algebra sobre (A,s) cuando A es un producto cruzado} is equivalent to items~(2)--(6). Item~(4) of Theorem~\ref{condiciones para la existencia de una estructura de H_q modulo algebra sobre (A,s) cuando A es un producto cruzado} and two of the equalities in item~(6) of the same theorem, are trivially satisfied because $\ov{\delta}_1 = \ov{\delta}_2 = 0$. Since
$$
\hat{\delta}_i\bigl(\hat{\alpha}(x_i)\bigr) = \nu_i\hat{\delta}_i(x_i) = \nu P_iw_{g_i} \quad \text{and}\quad \alpha\bigl(\hat{\delta}_i(x_i)\bigr) = \alpha(P_iw_{g_i}) = \alpha(P_i)\chi_{\alpha}(g_i)w_{g_i},
$$
item~(3) of Theorem~\ref{condiciones para la existencia de una estructura de H_q modulo algebra sobre (A,s) cuando A es un producto cruzado} is true if and only if item~(7) and the second equality in item~(9) hold. Since $\hat{\alpha}$ is $k[G]$-linear, item~(5) of Theorem~\ref{condiciones para la existencia de una estructura de H_q modulo algebra sobre (A,s) cuando A es un producto cruzado} is an immediate consequence of item~(2) of Theorem~\ref{condiciones para la existencia de estructuras de H_q-modulo algebra (caso 1)}. Finally we consider the non-trivial equalities in item~(6) of Theorem~\ref{condiciones para la existencia de una estructura de H_q modulo algebra sobre (A,s) cuando A es un producto cruzado}. It is easy to see that $\hat{\delta}_i\bigl(\hat{\varsigma}(x_i)\bigr) = q \varsigma\bigl(\hat{\delta}_i(x_i)\bigr)$ if and only if the first equality in item~(9) holds. On the other hand $\hat{\delta}_i\bigl( \hat{\varsigma}(v)\bigr) = q \varsigma\bigl(\hat{\delta}_i(v)\bigr)$ for all $v\in \ker\hat{\delta}_i$ if and only if $\hat{\varsigma}(\ker\hat{\delta}_i) \subseteq \ker\hat{\delta}_i$, which follows from items~(2), (4) and~(7). The equality $\delta_2\bigl( \hat{\delta}_1(v)\bigr) = \delta_1\bigl(\hat{\delta}_2(v)\bigr)$ is trivially satisfied for $v\in \ker\hat{\delta}_1\cap\ker\hat{\delta}_2$, and for $v\in \{x_1,x_2\}$ it is equivalent to item~(8). Lastly, the remaining equality coincides with item~(10).
\end{proof}

\begin{remark}\label{remark: algunas consecuencias} The following facts hold:

\begin{itemize}

\item[-] Since $\hat{\alpha}$ and $\hat{\varsigma}$ are bijective $k[G]$-linear maps, from item~(2) Theorem~\ref{condiciones para la existencia de estructuras de H_q-modulo algebra (caso 1)} it follows that
\begin{equation}
{}^{g_1^{-1}}\! v = {}^{g_2}\! v\quad\text{for all $v\in \ker\hat{\delta}_1 \cap
\ker\hat{\delta}_2$.}\label{eq3}
\end{equation}

\smallskip

\item[-] Since $x_1 \in \ker\hat{\delta}_2$ and $\ker\hat{\delta}_2$ is $G$-stable, ${}^g \! x_1 - \lambda_{1g}x_1\in \ker\hat{\delta}_1\cap \ker\hat{\delta}_2$. Similarly ${}^g \! x_2 - \lambda_{1g}x_2\in \ker\hat{\delta}_1\cap \ker\hat{\delta}_2$.

\smallskip

\item[-] Since $\ker \hat{\delta}_i$ is a $G$-submodule of $V$ and the $k$-linear map
$$
\xymatrix @R=-2pt {V\rto & V\\
v \ar@{|->}[0,1] & {}^g\!v}
$$
is an isomorphism for each $g\in G$, it is impossible that ${}^g\! x_i\in \ker \hat{\delta}_i$. Consequently, $\lambda_{ig}\in k^{\times}$ for each $g\in G$. Moreover, using again that $\ker \hat{\delta}_i$ is a $G$-submodule of $V$, it is easy to see that the map $g\mapsto\lambda_{ig}$ is a group homomorphism. Items~(1), (2), (4), (7) and the fact that $\hat{\alpha}$ is bijective imply that also $\omega_1,\omega_2,\nu_1,\nu_2\in k^{\times}$.

\smallskip

\item[-] Since
$$
\hat{\varsigma}(x_1) = \hat{\alpha} ({}^{g_2}\!x_1) \equiv \lambda_{1g_2} \hat{\alpha}(x_1) \equiv \lambda_{1g_2} \nu_1 x_1 \pmod{\ker \hat{\delta}_1},
$$
we have $\omega_1 = \lambda_{1g_2} \nu_1$. A similar argument shows that $\nu_2 = \lambda_{2g_1} \omega_2$.

\end{itemize}
\end{remark}

\begin{corollary}\label{determinante} Assume that the conditions above Theorem~\ref{condiciones para la existencia de estructuras de H_q-modulo algebra (caso 1)} are fulfilled and that there exists an $H_q$-module algebra structure on $(A,s)$ satisfying
$$
\sigma \!\cdot\!v = \hat{\varsigma}(v),\quad \sigma \!\cdot\!w_g=\chi_{\varsigma}(g)w_g, \quad D_i \!\cdot\!v = \hat{\delta}_i(v) \quad\text{and}\quad D_i \!\cdot\!w_g =0
$$
for all $v\in V$, $g\in G$ and $i\in \{1,2\}$. If $P_1\in S(\ker \hat{\delta}_1)$ and $P_2\in S(\ker \hat{\delta}_2)$, then
$$
\lambda_{1g_1}\lambda_{1g_2}=q \quad\text{and}\quad \lambda_{2g_1}\lambda_{2g_2}=q^{-1}.
$$
Moreover $g_0 := g_1g_2$ has determinant~$1$ as an operator on $V$.
\end{corollary}

\begin{proof} By items~(9), (2) and~(5) of Theorem~\ref{condiciones para la existencia de estructuras de H_q-modulo algebra (caso 1)},
$$
q^{-1}\omega_1\chi^{-1}_{\varsigma}(g_1)P_1 =\varsigma(P_1) = {}^{g_1^{-1}}\!\hat{\alpha}(P_1)= \nu_1\chi^{-1}_{\alpha}(g_1){}^{g_1^{-1}}\! P_1=\nu_1 \lambda_{1g_1}^{-1}\chi^{-1}_{\varsigma} (g_1) P_1.
$$
Hence $\lambda_{1g_1}\lambda_{1g_2}=q$ as we want, since $\omega_1=\nu_1\lambda_{1g_2}$. The proof that $\lambda_{2g_1}\lambda_{2g_2}=q^{-1}$ is similar. It remains to check that $\det(g_0)=1$. Since $\ker \hat{\delta}_1$ and $\ker \hat{\delta}_2$ are $G$-invariant, we have
$$
{}^g\!x_1\in \ker \hat{\delta}_2\quad\text{and}\quad {}^g\!x_2\in \ker \hat{\delta}_1\quad\text{for all $g\in G$,}
$$
and so
$$
{}^{g_0}\!x_1\in \lambda_{1g_1}\lambda_{1g_2}x_1 +W\quad\text{and}\quad {}^{g_0}\!x_2\in \lambda_{2g_1}\lambda_{2g_2}x_1 +W,
$$
where $W=\ker \hat{\delta}_1\cap \ker \hat{\delta}_2$. Moreover, by Remark~\ref{remark: algunas consecuencias} we know that $g_0$ acts as the identity map on $W$ and hence $\det (g_0) = \lambda_{1g_1}\lambda_{1g_2}\lambda_{2g_1}\lambda_{2g_2}=1$.
\end{proof}

\begin{remark}\label{caso particular} A particular case is the $H_q$-module algebra $A$ considered in~\cite[Section~4]{W}, in which $P_1=1$, $g_1=1$ and $\hat{\alpha}$ is the identity map. Our $P_2$, $g_2$ and $f$ correspond in~\cite{W} to $s$, $g$ and $\alpha$, respectively. Our computations show that the condition that  $h(s) = x_1(h)x_2(h)\alpha(g,h) \alpha^{-1}(h,g)s$, which appears as informed by the cohomology of finite groups in~\cite{W}, is in fact necessary for the existence of the $H_q$-module algebra structure of $A$, and it does not depend on cohomological considerations. In particular we need this condition for any group $G$, finite or not. Similarly the conditions that $g$ is central and $\det(g)=1$ are necessary even for infinite groups.
\end{remark}

Let $G$, $V$, $f\colon G\times G\to k^{\times}$ and $A$ be as at the beginning of this section. Let $\hat{\alpha}\colon V\to V$ be a bijective $k[G]$-linear map, $\chi_{\alpha}\colon G\to k^{\times}$ a group homomorphism, $\alpha\colon A\to A$ the algebra automorphism induced by $\hat{\alpha}$ and $\chi_{\alpha}$, and $s$ the good transposition associated with $\alpha$. Let

\begin{itemize}

\smallskip

\item[a)] $V_1\neq V_2$ subspaces of codimension $1$ of $V$ such that $V_1$ and $V_2$ are $\hat{\alpha}$-stable $G$-submodules of $V$,

\smallskip

\item[b)] $g_1$ and $g_2$ central elements of $G$ such that ${}^{g_1^{-1}}\! v = {}^{g_2}\!v$ for all $v\in V_1\cap V_2$,

\smallskip

\item[c)] $\chi_{\varsigma}\colon G\to k^{\times}$ a group homomorphism and $\hat{ \varsigma}\colon V\to V$ the map defined by
$$
\hat{\varsigma}(v) := \begin{cases} \hat{\alpha}({}^{g_1^{-1}}\! v) &\quad \text{if $v\in V_1$,}\\ \hat{\alpha}({}^{g_2}\!v) &\quad\text{if $v\in V_2$,}\end{cases}
$$

\smallskip

\item[d)] $x_1\in V_2\setminus V_1$, $x_2\in V_1\setminus V_2$, $P_1\in S(V_1)$, $P_2\in S(V_2)$  and $\hat{\delta}_1,\hat{\delta}_2\colon V\to A$ the maps defined by
$$
\ker \hat{\delta}_i := V_i\quad\text{and}\quad \hat{\delta}_i(x_i) := P_iw_{g_i}.
$$

\end{itemize}

For $g\in G$ and $i\in \{1,2\}$, let $\lambda_{ig},\nu_i,\omega_i\in k^{\times}$ be the elements defined by the conditions ${}^g\! x_i - \lambda_{ig} x_i\in V_i$, $\hat{\alpha}(x_i) - \nu_ix_i\in V_i$ and $\hat{\varsigma}(x_i) - \omega_ix_i\in V_i$.

\smallskip

The following result is a sort of a reformulation of Theorem~\ref{condiciones para la existencia de estructuras de H_q-modulo algebra (caso 1)}, more appropriate to construct explicit examples. The only new hypothesis that we need is that $P_i\in S(V_i)$.

\begin{corollary}\label{coro: first case con P_i en S(V_i)} There is an $H_q$-module algebra
structure on $(A,s)$, satisfying
$$
\sigma \!\cdot\!v = \hat{\varsigma}(v),\quad \sigma \!\cdot\!w_g=\chi_{\varsigma}(g)w_g, \quad D_i \!\cdot\!v = \hat{\delta}_i(v) \quad\text{and}\quad D_i \!\cdot\!w_g=0
$$
for all $v\in V$, $g\in G$ and $i\in \{1,2\}$, if and only if

\begin{enumerate}

\smallskip

\item $q = \lambda_{1g_1}\lambda_{1g_2}$ and $q^{-1} = \lambda_{2g_1}\lambda_{2g_2}$,

\smallskip

\item ${}^g\! P_1 = \lambda_{1g} \chi_{\alpha}^{-1}(g)\chi_{\varsigma}(g)f^{-1}(g,g_1) f(g_1,g) P_1$,

\smallskip

\item ${}^g\! P_2 = \lambda_{2g}\chi_{\alpha}(g)\chi_{\varsigma}^{-1}(g)f^{-1}(g,g_2) f(g_2,g)P_2$,

\smallskip

\item $\alpha(P_i)=\nu_i\chi_{\alpha}^{-1}(g_i)P_i$,

\smallskip

\item $P_1\in \ker\delta_2$ and $P_2\in \ker\delta_1$, where $\delta_1,\delta_2\colon A\to A$ are the maps defined in item~(8) of Theorem~\ref{condiciones para la existencia de estructuras de H_q-modulo algebra (caso 1)},

\smallskip

\item If $q\ne 1$ and $q^l=1$, then $\delta_1^l = \delta_2^l = 0$.

\end{enumerate}

\end{corollary}

\begin{proof} $\Leftarrow$) \enspace By~a), b), c) and d), it is obvious that items~(1), (2), (3), (4) and~(7) of Theorem~\ref{condiciones para la existencia de estructuras de H_q-modulo algebra (caso 1)} are satisfied. Moreover items~(2), (3), (5) and~(6) are items~(5), (6), (8) and~(10) of Theorem~\ref{condiciones para la existencia de estructuras de H_q-modulo algebra (caso 1)}.  So, we only must to check that item~(9) of Theorem~\ref{condiciones para la existencia de estructuras de H_q-modulo algebra (caso 1)} is satisfied. But the second equality in this item is exactly the one required in item~(4) of the present corollary, and we are going to check that the first one is true with $q = \lambda_{1g_1}\lambda_{1g_2}$. Arguing as in Remark~\ref{remark: algunas consecuencias}, and using item~(2) with $g = g_1$, item~(1) and item~(4), we obtain
\begin{align*}
q^{-1}\omega_1 \chi_{\varsigma}^{-1}(g_1)P_1 &= q^{-1}\lambda_{1g_2}\nu_1 \chi_{\varsigma}^{-1} (g_1)P_1\\
&= q^{-1}\lambda_{1g_1}\lambda_{1g_2}\nu_1 \chi_{\alpha}^{-1}(g_1)\,{}^{g^{-1}_1}\! P_1\\
& =\nu_1 \chi_{\alpha}^{-1}(g_1)\,{}^{g^{-1}_1}\!P_1\\
&={}^{g^{-1}_1}\!\alpha(P_1)\\
&=\varsigma(P_1),
\end{align*}
where the last equality is true since $P_1\in S(V_1)$. Again arguing as in Remark~\ref{remark: algunas consecuencias}, and using item~(3) with $g = g_2$, item~(1) and item~(4), we obtain
\begin{align*}
q^{-1}\omega_2 \chi_{\varsigma}^{-1}(g_2)P_2 & = q^{-1}\lambda_{2g_1}^{-1}\nu_2 \chi_{\varsigma}^{-1}(g_2)P_2\\
& = q^{-1}\lambda_{2g_1}^{-1}\lambda_{2g_2}^{-1}\nu_2\chi_{\alpha}^{-1}(g_2)\,{}^{g_2}\!P_2\\
& = \nu_2\chi_{\alpha}^{-1}(g_2)\,{}^{g_2}\!P_2\\
& ={}^{g_2}\!\alpha(P_2)\\
& =\varsigma(P_2),
\end{align*}
where the last equality is true since $P_2\in S(V_2)$.

\smallskip

\noindent $\Rightarrow$)\enspace Items~(2), (3), (5) and~(6) are items (5), (6), (8) and~(1) of Theorem~\ref{condiciones para la existencia de estructuras de H_q-modulo algebra (caso 1)}, and item~(4) is the first  equality in item~(9) of that theorem. Finally item~(1) follows from Corollary~\ref{determinante}.
\end{proof}

The following result shows that if $x_1$ and $x_2$ are eigenvectors of the maps $v\mapsto {}^{g_1}\! v$ and $v\mapsto {}^{g_2}\! v$, then item~(5) in the statement of Corollary~\ref{coro: first case con P_i en S(V_i)} can be easily tested and item~(6) can be removed from the hypothesis.

\begin{proposition}\label{variante del corolario} Assume that conditions~a), b), c) and~d) above Corollary~\ref{coro: first case con P_i en S(V_i)} are fulfilled. Let $\delta_1$ and $\delta_2$ be the maps introduced in item~(8) of Theorem~\ref{condiciones para la existencia de estructuras de H_q-modulo algebra (caso 1)}. If
$$
\lambda_{1g_1}\lambda_{1g_2}=q,\quad \lambda_{2g_1}\lambda_{2g_2}=q^{-1}\quad\text{and}\quad {}^{g_i}\!x_j=\lambda_{jg_i}x_j \quad\text{for $1\le i,j\le 2$,}
$$
then:
\begin{enumerate}

\smallskip

\item $\delta_1^l = \delta_2^l = 0$, whenever $q\ne 1$ and $q^l=1$.

\smallskip

\item If $q=1$ or it is not a root of unity, then $P_1\in \ker \delta_2$ and $P_2\in \ker \delta_1$ if and only if $P_1,P_2\in S(V_1\cap V_2)$.

\smallskip

\item If $q\ne 1$ is a primitive $l$-root of unity, then $P_1\in \ker \delta_2$ and $P_2\in \ker \delta_1$ if and only if $P_1\in S\bigl(k\, x_2^l \oplus (V_1\cap V_2)\bigr)$ and $P_2\in S\bigl(k\, x_1^l \oplus (V_1\cap V_2)\bigr)$.

\end{enumerate}

\end{proposition}

\begin{proof}  The proposition is a direct consequence of the following formulas:
\begin{align*}
& \delta_1^s(x_1^{r_1}\cdots x_n^{r_n}w_g) = \begin{cases} c \alpha^s\bigl(x_1^{r_1-s} x_2^{r_2}\cdots x_n^{r_n}\bigr)w_{g_1^sg} &\text{for $s\le r_1$,}\\ 0 &\text{otherwise,} \end{cases}
\intertext{and}
& \delta_2^s(x_1^{r_1}\cdots x_n^{r_n}w_g) = \begin{cases}  d x_2^{r_2-s}{}^{g_2^s}\! \bigl(x_1^{r_1}x_3^{r_3} \cdots x_n^{r_n}\bigr) w_{g_2^sg} &\text{for $s\le r_2$,} \\ 0 & \text{otherwise,}\end{cases}
\end{align*}
where $\alpha^s$ denotes the $s$-fold composition of $\alpha$,
\begin{align*}
& c=\chi_{\varsigma}^s(g)\chi_{\varsigma}^{s(s-1)/2}(g_1)\chi_{\alpha}^{s(s-1)/2}(g_1) \Biggl( \prod_{k=0}^{s-1} (r_1-k)_q\Biggr)\Biggl(\prod_{k=0}^{s-1} f(g_1,g_1^kg)\Biggr) \alpha^{s-1} (P_1^s),\\
& d=\lambda_{2g_2}^{sr_2-s(s+1)/2}\Biggl(\prod_{k=0}^{s-1} (r_2-k)_q\Biggr)\Biggl( \prod_{k=0}^{s-1} f(g_2,g_2^kg)\Biggr) \Biggl(\prod_{k=0}^{s-1} {}^{g_2^k}\!P_2\Biggr).
\end{align*}
We will prove the formula for $\delta_1^s$ and we will leave the other one to the reader. We begin with the case $s=1$. Since $x_2,\dots,x_n\in \ker \hat{\delta}_1$ and $\hat{\delta}_1 (x_1) = P_1w_{g_1}$, from the definition of $\delta_1$ it follows that
$$
\delta_1(x_1^{r_1}\cdots x_n^{r_n}w_g) = \sum_{j=0}^{r_1-1} \alpha(x_1^j) P_1w_{g_1} \varsigma (x_1^{r_1-j-1}x_2^{r_2}\cdots x_n^{r_n}w_g).
$$
Thus, using the definition of $\varsigma$,  item~c) above Corollary~\ref{coro: first case con P_i en S(V_i)}, the fact that $\alpha$ is $G$-linear and the hypothesis, we obtain
\begin{align*}
\delta_1(x_1^{r_1}\cdots x_n^{r_n}w_g) & = \sum_{j=0}^{r_1-1} \alpha(x_1^j) P_1w_{g_1} \alpha ({}^{g_2}\! x_1^{r_1-j-1}) {}^{g_1^{-1}}\!\alpha (x_2^{r_2}\cdots x_n^{r_n}) \chi_{\varsigma}(g)w_g\\
& = \sum_{j=0}^{r_1-1} \alpha(x_1^j) P_1 \alpha ({}^{g_1g_2}\! x_1^{r_1-j-1}) \alpha (x_2^{r_2}\cdots x_n^{r_n}) \chi_{\varsigma}(g)f(g_1,g) w_{g_1g}\\
& = \chi_{\varsigma}(g) (r_1)_q f(g_1,g) P_1  \alpha (x_1^{r_1-1}x_2^{r_2}\cdots x_n^{r_n}) w_{g_1g}.
\end{align*}
Assume that $s\le r_1$ and that the formula for $\delta_1^s$ holds. Since $c$ depends on $s$, $r_1$ and $g$, it will be convenient for us to use  the more precise notation  $c_{s,r_1}(g)$ for $c$. From items~(3) and~(5) of Theorem~\ref{condiciones para la existencia de una estructura de H_q modulo algebra sobre (A,s) cuando A es un producto cruzado} and item~(9) of Theorem~\ref{algunas acciones de H_q sobre algebras}, It follows easily that $\alpha\xcirc\delta_1 = \delta_1\xcirc \alpha$ on $S(V)$. Using this fact, item~(9) of Theorem~\ref{algunas acciones de H_q sobre algebras} and the inductive hypothesis, we obtain
$$
\delta_1^{s+1}(x_1^{r_1}\cdots x_n^{r_n}w_g) = \alpha\bigl(c_{sr_1}(g)\bigr)\alpha^s\bigl( \delta_1 (x_1^{r_1-s} x_2^{r_2}\cdots x_n^{r_n})\bigr) \varsigma(w_{g_1^sg}).
$$
If $s=r_1$, then $\delta_1 (x_1^{r_1-s} x_2^{r_2}\cdots x_n^{r_n})=0$. Otherwise,
\begin{align*}
\delta_1^{s+1}(x_1^{r_1}\cdots x_n^{r_n}w_g) & = \ov{c}\alpha^s\bigl(\alpha(x_1^{r_1-s-1} x_2^{r_2}\cdots x_n^{r_n})w_{g_1}\bigr) \varsigma(w_{g_1^sg})\\
& = \ov{c}\alpha^{s+1}(x_1^{r_1-s-1} x_2^{r_2}\cdots x_n^{r_n})\chi_{\alpha}^s(g_1) \chi_{\varsigma}^s(g_1)\chi_{\varsigma}(g) f(g_1,g_1^sg) w_{g_1^{s+1}g},
\end{align*}
where $\ov{c} = \alpha\bigl(c_{s,r_1}(g)\bigr)\alpha^s\bigl(c_{1,r_1-s}(1) \bigr)$. The formula for $\delta_1^{s+1}$ follows immediately from this fact.
\end{proof}

\begin{example}\label{estructuras de módulo algebra con G cíclico} Let $G=\langle g
\rangle$ be an order $r$ cyclic group, $\xi$ an element of $k^{\times}$ and $f_{\xi}\colon G\otimes G \to k$ the cocycle defined by
$$
f_{\xi}(g^u,g^v) := \begin{cases} 1& \text{if $u+v<r$,}\\ \xi &\text{otherwise.}\end{cases}
$$
Of course, if $r=\infty$, then for any $\xi$ this is the trivial cocycle. Let $V$ be a vector space endowed with an action of $G$ and let $A$ be the crossed product $A=S(V)\#_{f_{\xi}} G$. Let $\{x_1,\dots,x_n\}$ be a basis of $V$. Let us $V_1$ and $V_2$ denote the subspaces of $V$ generated by $\{x_2,\dots,x_n\}$ and $\{x_1,x_3,\dots,x_n\}$ respectively. Let $\hat{\alpha}\colon V\to V$ be a bijective $k[G]$-linear map. Assume that $V_1$ and $V_2$ are $\hat{\alpha}$-stable $G$-submodules of $V$ and that there exist $\lambda_1,\lambda_2\in k^{\times}$ such that ${}^g\!x_1=\lambda_1x_1$ and ${}^g\!x_2=\lambda_2x_2$. Let $m_1,m_2\in \mathds{Z}$. Assume that ${}^{g^{m_1+m_2}}\!v = v$ for all $v\in V_1\cap V_2$ (if $r<\infty$ we can take $0\le m_1,m_2<r$). Let $\hat{\varsigma}\colon V\to V$ be the map defined by
$$
\hat{\varsigma}(v) :=\begin{cases}\hat{\alpha}({}^{g^{-m_1}}\! v) &\quad \text{if $v\in V_1$,}\\ \hat{\alpha}({}^{g^{m_2}}\!v) &\quad\text{if $v\in V_2$,} \end{cases}
$$
and let $\chi_{\alpha}, \chi_{\varsigma}\colon G\to k^{\times}$ be two morphisms. Consider the automorphism of algebras $\alpha\colon A\to A$ given by $\alpha(v):= \hat{\alpha}(v)$ for $v\in V$ and $\alpha(w_g)= \chi_{\alpha}(g)w_g$, and define $\hat{\delta}_1,\hat{\delta}_2\colon V\to A$ by
\begin{alignat*}{2}
&\hat{\delta}_1(x_2)=\cdots=\hat{\delta}_1(x_n):=0,&&\quad\hat{\delta}_1(x_1):=P_1w_{g^{m_1}},\\
&\hat{\delta}_2(x_1)=\hat{\delta}_2(x_3)=\cdots=\hat{\delta}_1(x_n):=0,&&\quad \hat{\delta}_2(x_2):=P_2w_{g^{m_2}},
\end{alignat*}
where $P_1\in S(V_1)\setminus\{0\}$ and $P_2\in S(V_2)\setminus\{0\}$. Let $s$ be the transposition of $H_q$ with $A$ associated with $\alpha$. There is an $H_q$-module algebra structure over $(A,s)$ satisfying
$$
\sigma\!\cdot\!v = \hat{\varsigma}(v),\quad\!\sigma\!\cdot\!w_g = \chi_{\varsigma}(g)w_g,\quad\! D_i\!\cdot\!v = \hat{\delta}_i(v)\quad\!\text{and}\quad\! D_i\!\cdot\!w_g = 0\quad\!\text{for all $v\in V$,}
$$
if and only if

\begin{enumerate}

\smallskip

\item $q=\lambda_1^{m_1+m_2}$ and $q^{-1}=\lambda_2^{m_1+m_2}$,

\smallskip

\item ${}^g\!P_1=\lambda_1\chi^{-1}_{\alpha}(g)\chi_{\varsigma}(g)P_1$ and ${}^g\!P_2 =\lambda_2\chi_{\alpha}(g)\chi_{\varsigma}^{-1}(g)P_2$,

\smallskip

\item $\alpha(P_1)=\nu_1\chi^{-m_1}_{\alpha}(g)P_1$ and $\alpha(P_2) = \nu_2 \chi^{-m_2}_{\alpha}(g) P_2$,

\smallskip

\item If $q = 1$ or $q$ is not a root of unity, then $P_1,P_2\in k[x_3,\dots,x_n]$,

\smallskip

\item If $q\ne 1$ is a primitive $l$-root of unity, then
$$
P_1\in k[x_2^l,x_3,\dots,x_n]\quad and\quad P_2\in k[x_1^l,x_3,\dots,x_n].
$$

\smallskip

\end{enumerate}
Consequently, in order to obtain explicit examples of braided $H_q$-module algebra structures on an algebra $A$ of the shape  $S(V)\#_{f_{\xi}} G$, where $V$ is a $k$-vector space with basis $\{x_1,\dots, x_n\}$ and $G=\langle g \rangle$ is a cyclic group of order $r\le \infty$, we proceed as follows:

\begin{description}

\smallskip

\item[First] We define an action of $G$ on $V$. For this we choose

\begin{itemize}

\smallskip

\item[-] a $k$-linear automorphism $\gamma$ of $V_{12}:=\langle x_3,\dots,x_n\rangle$, whose order divides $r$ if $r<\infty$,

\smallskip

\item[-] $\lambda_1,\lambda_2\in k^{\times}$ such that $\lambda_1^r = \lambda_2^r =1$ if $r<\infty$,

\smallskip

\end{itemize}
and we set
$$
{}^g\! x_i := \begin{cases} \lambda_1 x_1 & \text{if $i = 1$,}\\ \lambda_2 x_2 & \text{if $i = 2$,}\\ \gamma(x_i) & \text{if $i\ge 3$.}\end{cases}
$$

\smallskip

\item[Second] We construct the algebra $A$. For this we choose $\xi\in k^{\times}$ and we define $A = S(V)\#_{f_{\xi}} G$, where $f_{\xi}$ is the cocycle associate with $\xi$.

\smallskip

\item[Third] We endow $A$ with a $k$-algebra automorphism $\alpha$. For this we take $\nu_1,\nu_2, \eta\in k^{\times}$ such that $\eta^r=1$ if $r<\infty$, a $k$-linear automorphism $\alpha'$ of $V_{12}$ and $v_1,v_2\in V_{12}$, and we define
$$
\alpha(w_g) := \eta w_g\quad\text{and}\quad\alpha(x_i) := \begin{cases} \nu_1 x_1 + v_1& \text{if $i = 1$,}\\ \nu_2x_2 + v_2& \text{if $i = 2$,}\\ \alpha'(x_i) & \text{if $i\ge 3$.}\end{cases}
$$

\smallskip

\item[Fourth] We choose $m_1,m_2\in \mathds{Z}$ and $\zeta\in k^{\times}$ such that
$$
\gamma^{m_1+m_2} = \ide,\quad (\lambda_1\lambda_2)^{m_1+m_2} =1\quad\text{and}\quad \zeta^r=1\,\text{ if $r<\infty$,}
$$
and we define
$$
\varsigma(w_g) := \zeta w_g\quad\text{and}\quad\varsigma(x_i) := \begin{cases} \lambda_1^{m_2}(\nu_1 x_1 + v_1)& \text{if $i = 1$,}\\ \lambda_2^{-m_1}(\nu_2x_2 + v_2)& \text{if $i = 2$,}\\ \alpha'\bigl(\gamma^{m_2}(x_i)\bigr) & \text{if $i\ge 3$.}\end{cases}
$$

\smallskip

\item[Fifth] we set $q:=\lambda_1^{m_1+m_2}$ and we choose $P_1,P_2\in S(V)\setminus\{0\}$ such that

\begin{itemize}

\smallskip

\item[-] if $q$ is not a root of unity, then $P_1,P_2\in k[x_3,\dots,x_n]$,

\smallskip

\item[-] if $q$ is a primitive $l$-root of unity, then
$$
P_1\in k[x_2^l,x_3,\dots,x_n]\quad\text{and}\quad P_2\in k[x_1^l,x_3,\dots,x_n],
$$

\smallskip

\item[-] ${}^g\!P_1=\lambda_1\eta^{-1}\zeta P_1$ and ${}^g\!P_2 = \lambda_2\eta \zeta^{-1}P_2$,

\smallskip

\item[-] $\alpha(P_1)=\nu_1\eta^{-m_1}P_1$ and $\alpha(P_2) = \nu_2\eta^{-m_2}P_2$.

\end{itemize}

\end{description}
Now, by the discussion at the beginning of this example, there is an $H_q$-module algebra structure on $(A,s)$, where $s\colon H_q\ot A\to A\ot H_q$ is the good transposition associated with $\alpha$, such that
$$
\sigma\!\cdot\!x_j = \varsigma(x_j),\quad\!\sigma\!\cdot\!w_g = \zeta w_g,\quad\! D_i\!\cdot\!w_g = 0\quad \text{and}\quad D_i(x_j)=\begin{cases} 0 &\text{if $i\ne j$,}\\ P_i w_{g^{mi}} &\text{if $i= j$,}\end{cases}
$$
where $i\in\{1,2\}$ and $j\in \{1,\dots,n\}$.
\end{example}

\begin{remark} If $P_1(0)\ne 0$ and $P_2(0)\ne 0$, then the conditions in the first step are fulfilled if and only if $\lambda_1\lambda_2 = 1$, $\eta = \lambda_1\zeta$, $\nu_1 =\eta^{m_1}$, $\nu_2 =\eta^{m_2}$, $P_1$ and $P_2$ are $G$-invariants, $\alpha(P_1) = P_1$ and $\alpha(P_2) = P_2$.
\end{remark}

\subsection{Second case.}\label{Second case} Let $\hat{\alpha}$, $\chi_{\alpha}$, $\alpha$ and $s$ be as in the discussion following Proposition~\ref{prop: algunas transposiciones buenas}, let $\chi_{\varsigma}\colon G\to k^{\times}$ be a map and let $\hat{\delta}_1\colon V\to A$, $\hat{\delta}_2\colon V\to A$ and $\hat{\varsigma}\colon V\to V$ be $k$-linear maps such that $\ker \hat{\delta}_1\neq \ker \hat{\delta}_2$ are subspaces of codimension~$1$ of $V$. Here we are going to consider a more general situation that the one studied in the previous subsection. Assume that for each $i\in \{1,2\}$ there exist

\begin{itemize}

\smallskip

\item[-] an element $x_i\in V\setminus\ker(\hat{\delta}_i)$,

\smallskip

\item[-] different elements $g_{i1}, \dots, g_{in_i}$ of $G$,

\smallskip

\item[-] polynomials $P^{(i)}_{g_{i1}}, \dots, P^{(i)}_{g_{in_i}}\in S(V)\setminus \{0\}$,

\smallskip

\end{itemize}
such that
$$
\hat{\delta}_i(x_i) = \sum_{j=1}^{n_i} P^{(i)}_{g_{ij}}w_{g_{ij}}.
$$
(The reason for the notation $P^{(i)}_{g_{ij}}$ instead of the more simple $P_{ij}$ will became clear in items~(5) and (6) of the following theorem). Without loss of generality we can assume that $x_1\in \ker\hat{\delta}_2$ and $x_2\in \ker\hat{\delta}_1$ (and we do it). For $g\in G$ and $i\in\{1,2\}$, let $\lambda_{ig},\omega_i,\nu_i\in k$  be the elements defined by the following conditions:
$$
{}^g\! x_i - \lambda_{ig} x_i\in \ker \hat{\delta}_i,\quad \hat{\varsigma}(x_i) - \omega_i x_i\in \ker \hat{\delta}_i\quad\text{and}\quad \hat{\alpha}(x_i) - \nu_i x_i\in \ker\hat{\delta}_i.
$$

\begin{lemma}\label{buena def de delta_i caso 2} The following facts hold:

\begin{enumerate}

\smallskip

\item Condition~\eqref{eq2.1} of Lemma~\ref{buena def de delta_i} is satisfied if and only if ${}^{g_{1j}}\!\hat{\varsigma}(v) = \hat{\alpha}(v)$ for all $j\le n_1$ and $v\in \ker \hat{\delta}_1$.

\smallskip

\item Condition~\eqref{eq2.2} of Lemma~\ref{buena def de delta_i} is satisfied if and only if ${}^{g_{2j}}\! v = \hat{\varsigma}\bigl(\hat{\alpha}^{-1}(v) \bigr)$ for all $j\le n_2$ and $v\in \ker \hat{\delta}_2$.

\end{enumerate}

\end{lemma}

\begin{proof} Mimic the proof of Lemma~\ref{buena def de delta_i caso 1}.
\end{proof}

\begin{lemma}\label{propiedades de delta_i caso 2} The following facts hold:

\begin{enumerate}

\smallskip

\item Items~(1)(a) and~(1)(b) of Lemma~\ref{propiedades de delta_i} are satisfied if and only if

\begin{enumerate}

\smallskip

\item $\ker \hat{\delta}_1$ is a $G$-submodule of $V$,

\smallskip

\item $\{g_{1j}: 1\le j\le n_1\}$ is a union of conjugacy classes of $G$,

\smallskip

\item ${}^g\! P^{(1)}_{g_{1j}} = \lambda_{1g}\chi_{\alpha}^{-1}(g)\chi_{\varsigma}(g) f^{-1}(g,g_{1j}) f(gg_{1j}g^{-1},g) P^{(1)}_{gg_{1j}g^{-1}}$ for $j\le n_1$.

\smallskip

\end{enumerate}

\smallskip

\item Items~(2)(a) and~(2)(b) of Lemma~\ref{propiedades de delta_i} are satisfied if and only if

\begin{enumerate}

\smallskip

\item $\ker \hat{\delta}_2$ is a $G$-submodule of $V$,

\smallskip

\item $\{g_{2j}: 1\le j\le n_2\}$ is a union of conjugacy classes of $G$,

\smallskip

\item ${}^g\!P^{(2)}_{g_{2j}} = \lambda_{2g}\chi_{\alpha}(g)\chi_{\varsigma}^{-1}(g) f^{-1}(g,g_{2j})f(gg_{2j}g^{-1},g)P^{(2)}_{gg_{2j}g^{-1}}$ for $j\le n_2$.

\smallskip

\end{enumerate}
\end{enumerate}
\end{lemma}

\begin{proof} Mimic the proof of Lemma~\ref{propiedades de delta_i caso 1}.
\end{proof}

\begin{theorem}\label{condiciones para la existencia de estructuras de H_q-modulo algebra (caso 2)} There is an $H_q$-module algebra structure on $(A,s)$, satisfying
$$
\sigma\!\cdot\!v = \hat{\varsigma}(v),\quad \sigma \!\cdot\!w_g=\chi_{\varsigma}(g)w_g, \quad D_i \!\cdot\!v = \hat{\delta}_i(v) \quad\text{and}\quad D_i \!\cdot\!w_g =0
$$
for all $v\in V$, $g\in G$ and $i\in \{1,2\}$,  if and only if

\begin{enumerate}

\smallskip

\item $\hat{\varsigma}$ is a bijective $k[G]$-linear map and $\chi_{\varsigma}$ is a group homomorphism,

\smallskip

\item $\hat{\varsigma}(v) = {}^{g_{1j}^{-1}}\!\hat{\alpha}(v)$ for $j\le n_1$ and all $v\in \ker \hat{\delta}_1$, and  $\hat{\varsigma}(v) = {}^{g_{2j}}\!\hat{\alpha}(v)$ for $j\le n_2$ and all $v\in \ker \hat{\delta}_2$,

\smallskip

\item $\{g_{ij}: 1\le j\le n_i\}$ is a union of conjugacy classes of~$G$ for $i\in \{1,2\}$,

\smallskip

\item $\ker \hat{\delta}_1$ and $\ker \hat{\delta}_2$ are $G$-submodules of $V$,

\smallskip

\item ${}^g\! P^{(1)}_{g_{1j}} =\lambda_{1g}\chi_{\alpha}^{-1}(g)\chi_{\varsigma}(g) f^{-1}(g,g_{1j}) f(gg_{1j}g^{-1},g) P^{(1)}_{gg_{1j}g^{-1}}$ for $j\le n_1$,

\smallskip

\item  ${}^g\!P^{(2)}_{g_{2j}}=\lambda_{2g}\chi_{\alpha}(g)\chi_{\varsigma}^{-1}(g) f^{-1}(g,g_{2j}) f(gg_{2j}g^{-1},g)P^{(2)}_{gg_{2j}g^{-1}}$ for $j\le n_2$,

\smallskip

\item $\hat{\alpha}(\ker \hat{\delta}_i) = \ker \hat{\delta}_i$ for $i\in \{1,2\}$,

\smallskip

\item $\sum_{j=1}^{n_1}P^{(1)}_{g_{1j}}w_{g_{1j}}\in \ker\delta_2$ and $\sum_{j=1}^{n_2} P^{(2)}_{g_{2j}} w_{g_{2j}}\in \ker\delta_1$, where $\delta_1$ and $\delta_2$ are the maps defined by
\begin{align*}
&\qquad\delta_1(\mathbf{v}_{1m} w_g) := \sum_{j=1}^m \alpha(\mathbf{v}_{1,j-1}) \hat{\delta}_1(v_j) \varsigma(\mathbf{v}_{j+1,m}w_g),\\
&\qquad\delta_2(\mathbf{v}_{1m} w_g) := \sum_{j=1}^m\varsigma\bigl(\alpha^{-1} (\mathbf{v}_{1,j-1})\bigr) \hat{\delta}_2(v_j)\mathbf{v}_{j+1,m}w_g,
\end{align*}
in which $\mathbf{v}_{hl} = v_h\cdots v_l$,

\smallskip

\item $\varsigma(P^{(i)}_{g_{ij}}) = q^{-1} \omega_i\chi_{\varsigma}^{-1} (g_{ij}) P^{(i)}_{g_{ij}}$ and $\alpha(P^{(i)}_{g_{ij}}) = \nu_i\chi_{\alpha}^{-1}(g_{ij}) P^{(i)}_{g_{ij}}$ for $i\in \{1,2\}$ and $j\le n_i$, where  $\varsigma$ is the map
    given by
$$
\varsigma(\mathbf{v}_{1m}w_g) := \hat{\varsigma}(v_1)\cdots\hat{\varsigma}(v_m) \chi_{\varsigma}(g) w_g.
$$

\smallskip

\item If $q\ne 1$ and $q^l=1$, then $\delta_1^l=\delta_2^l=0$.

\smallskip

\end{enumerate}

\end{theorem}

\begin{proof} Mimic the proof of Theorem~\ref{condiciones para la existencia de estructuras de H_q-modulo algebra (caso 1)}, but using Lemmas~\ref{buena def de delta_i caso 2} and~\ref{propiedades de delta_i caso 2} instead of Lemmas~\ref{buena def de delta_i caso 1} and~\ref{propiedades de delta_i caso 1}, respectively.
\end{proof}

\begin{remark}\label{remark: algunas consecuencias 2} Since $\alpha$ and $\varsigma$ are bijective $k[G]$-linear maps, from item~(2) it follows that
\begin{xalignat}{2}
& {}^{g_{1j}}\! v={}^{g_{1h}}\! v &&\text{for $1\le j,h\le n_1$ and all $v\in \ker\hat{\delta}_1$,}\label{eq6}\\
& {}^{g_{2j}}\! v={}^{g_{2h}}\! v &&\text{for $1\le j,h\le n_2$ and all $v\in \ker\hat{\delta}_2$,}\label{eq7}\\
&{}^{g_{1j}^{-1}}\! v = {}^{g_{2h}}\! v &&\text{for $1\le j\le n_1$, $1\le h\le n_2$ and all $v\in \ker\hat{\delta}_1 \cap \ker\hat{\delta}_2$.}\label{eq8}
\end{xalignat}
On the other hand, arguing as in Remark~\ref{remark: algunas consecuencias} we can check that

\begin{itemize}

\smallskip

\item[-] ${}^g \! x_i - \lambda_{1g}x_i\in\ker\hat{\delta}_1\cap \ker\hat{\delta}_2$ for all $g\in G$,

\smallskip

\item[-] $\lambda_{ig}\in k^{\times}$ for all $g\in G$,

\smallskip

\item[-] the maps $g\mapsto\lambda_{ig}$ are morphisms,

\smallskip

\item[-]  $\omega_1,\omega_2,\nu_1,\nu_2\in k^{\times}$.

\smallskip

\end{itemize}
Finally, since
$$
\hat{\varsigma}(x_1) = \hat{\alpha}({}^{g_{2j}}\!x_1)\equiv\lambda_{1g_{2j}} \hat{\alpha}(x_1) \pmod{\ker\hat{\delta}_1},
$$
we have $\omega_1 = \lambda_{1g_{2j}}\nu_1$ for $j\le n_2$. Similarly, $\nu_2 = \lambda_{2g_{1j}} \omega_2$ for $j \le n_1$. Consequently,
$$
\lambda_{1g_{21}} = \cdots = \lambda_{1g_{2n_2}}\quad\text{and}\quad \lambda_{2g_{11}}=\cdots = \lambda_{2g_{1n_1}},
$$
which also follows from~\eqref{eq6} and~\eqref{eq7}.
\end{remark}

\begin{corollary}\label{determinante (caso 2)} Assume that the conditions at the beginning of the present subsection are fulfilled and that there exists an $H_q$-module algebra structure on $(A,s)$, satisfying
$$
\sigma \!\cdot\!v = \hat{\varsigma}(v),\quad \sigma \!\cdot\!w_g = \chi_{\varsigma}(g)w_g, \quad D_i \!\cdot\!v = \hat{\delta}_i(v) \quad\text{and}\quad D_i \!\cdot\!w_g = 0
$$
for all $v\in V$, $g\in G$ and $i\in \{1,2\}$. If $P^{(1)}_{g_{1j}}\in S(\ker \hat{\delta}_1)$ and $P^{(2)}_{g_{2h}}\in S(\ker \hat{\delta}_2)$ for all $j\le n_1$ and $h\le n_2$, then
$$
\lambda_{1g_{1j}}\lambda_{1g_{2h}}=q\quad\text{and}\quad\lambda_{2g_{1j}} \lambda_{2g_{2h}}= q^{-1}.
$$
Moreover $g_{1j}g_{2h}$ has determinant~$1$ as an operator on $V$.
\end{corollary}

\begin{proof} This result generalizes Corollary~\ref{determinante}, and its proof is similar.
\end{proof}

Let $G$, $V$, $f\colon G\times G\to k^{\times}$, $A$, $\hat{\alpha}\colon V\to V$, $\chi_{\alpha}\colon G\to k^{\times}$, $\alpha\colon A\to A$ and $s$ be as below of Remark~\ref{caso particular}. Assume we have

\begin{itemize}

\smallskip

\item[a)] subspaces $V_1\neq V_2$ of codimension $1$ of $V$ such that $V_1$ and $V_2$ are $\hat{\alpha}$-stable $G$-submodules of $V$, and vectors $x_1\in V_2\setminus V_1$ and $x_2\in V_1\setminus V_2$,

\smallskip

\item[b)] different elements $g_{i1},\dots,g_{in_i}$ of $G$, where $i\in\{1,2\}$, such that:

\begin{itemize}

\smallskip

\item $\{g_{11},\dots,g_{1n_1}\}$ and $\{g_{21},\dots g_{2n_2}\}$ are unions of conjugacy classes of $G$,

\smallskip

\item ${}^{g_{1j}}\! v={}^{g_{1h}}\! v$ for $1\le j,h\le n_1$ and all $v\in V_1$,

\smallskip

\item ${}^{g_{2j}}\! v={}^{g_{2h}}\! v$ for $1\le j,h\le n_2$ and all $v\in V_2$,

\smallskip

\item ${}^{g_{1j}^{-1}}\! v = {}^{g_{2h}}\! v$ for $1\le j\le n_1$, $1\le h\le n_2$ and all $v\in V_1 \cap V_2$.
\end{itemize}

\smallskip

\item[c)] a morphism $\chi_{\varsigma}\colon G\to k^{\times}$

\smallskip

\item[d)]  Nonzero polynomials $P^{(1)}_{g_{1j}}\in S(V_1)$ and $P^{(2)}_{g_{2h}}\in S(V_2)$, where $1\le j\le n_1$ and $1\le h\le n_2$.

\smallskip

\end{itemize}
Let $\hat{\varsigma}\colon V\to V$ and $\hat{\delta}_1,\hat{\delta}_2 \colon V\to A$ be the maps defined by
$$
\hat{\varsigma}(v):=\begin{cases} \hat{\alpha}({}^{g_{11}^{-1}}\! v) &\quad \text{if $v\in V_1$,}\\\hat{\alpha}({}^{g_{21}}\!v) &\quad\text{if $v\in V_2$,} \end{cases} \quad \ker \hat{\delta}_i := V_i\quad\text{and}\quad \hat{\delta}_i(x_i) := \sum_{j=1}^{n_i} P^{(i)}_{g_{ij}} w_{g_{ij}}.
$$
For $g\in G$ and $i\in \{1,2\}$, let $\lambda_{ig},\nu_i\in k^{\times}$ be the elements defined by the following conditions: ${}^g\! x_i - \lambda_{ig} x_i\in V_i$ and $\alpha(x_i) - \nu_ix_i\in V_i$. Note that, by item~b),
$$
\lambda_{2g_{11}}=\cdots = \lambda_{2g_{1n_1}}\quad\text{and} \quad\lambda_{1g_{21}}=\cdots =\lambda_{1g_{2n_2}}.
$$

\begin{corollary}\label{coro: second case con P_i en S(V_i)} There is an $H_q$-module algebra structure on $(A,s)$, satisfying
$$
\sigma \!\cdot\!v = \hat{\varsigma}(v),\quad \sigma \!\cdot\!w_g=\chi_{\varsigma}(g)w_g, \quad D_h \!\cdot\!v = \hat{\delta}_h(v) \quad\text{and}\quad D_h \!\cdot\!w_g=0,
$$
for all $v\in V$, $g\in G$ and $i\in \{1,2\}$, if and only if for all $j\le n_1$ and $h\le n_2$ the following facts hold:

\begin{enumerate}

\smallskip

\item $q = \lambda_{1g_{1j}}\lambda_{1g_{21}}$ and $q^{-1} = \lambda_{2g_{11}} \lambda_{2g_{2h}}$,

\smallskip

\item ${}^g\! P^{(1)}_{g_{1j}} = \lambda_{1g}\chi_{\alpha}^{-1}(g)\chi_{\varsigma}(g) f^{-1}(g,g_{1j}) f(gg_{1j}g^{-1},g) P^{(1)}_{gg_{1j}g^{-1}}$,

\smallskip

\item ${}^g\! P^{(2)}_{g_{2h}} = \lambda_{2g}\chi_{\alpha}(g)\chi_{\varsigma}^{-1}(g) f^{-1}(g,g_{2h}) f(gg_{2h}g^{-1},g)P^{(2)}_{gg_{2h}g^{-1}}$,

\smallskip

\item $\alpha(P^{(1)}_{g_{1j}})=\nu_1\chi_{\alpha}^{-1}(g_{1j})P^{(1)}_{g_{1j}}$ and $\alpha(P^{(2)}_{g_{2h}})=\nu_2\chi_{\alpha}^{-1}(g_{2h})P^{(2)}_{g_{2h}}$,

\smallskip

\item $\sum_{j=1}^{n_1} P^{(1)}_{g_{1j}}w_{g_{1j}}\in \ker\delta_2$ and $\sum_{h=1}^{n_2} P^{(2)}_{g_{2h}}w_{g_{2h}}\in \ker \delta_1$, where $\delta_1,\delta_2\colon A\to A$ are the maps defined in item~(8) of Theorem~\ref{condiciones para la existencia de estructuras de H_q-modulo algebra (caso 2)}

\smallskip

\item If $q\ne 1$ and $q^l=1$, then $\delta_1^l = \delta_2^l = 0$.

\end{enumerate}

\end{corollary}

\begin{proof} It is similar to the proof of Corollary~\ref{coro: first case con P_i en S(V_i)}, using Theorem~\ref{condiciones para la existencia de estructuras de H_q-modulo algebra (caso 2)} instead of Theorem~\ref{condiciones para la existencia de estructuras de H_q-modulo algebra (caso 1)}. The proof that $\varsigma$ is $G$-linear requires additionally the fact that ${}^{gg_{ij}g^{-1}}\!v= {}^{g_{ij}}\!v$ for $1\le i\le 2$ and $1\le j\le n_i$, which is true by b).
\end{proof}

\begin{remark}\label{remark: igualdad de los lambda_ij cuando j se mueve} Assume that the hypothesis of Corollary~\ref{coro: second case con P_i en S(V_i)} are fulfilled. Then, as it was note above this corollary,
$$
\lambda_{2g_{11}}=\cdots = \lambda_{2g_{1n_1}}\quad\text{and}\quad \lambda_{1g_{21}}=\cdots = \lambda_{1g_{2n_2}}.
$$
Moreover, by item~(1) it is clear that
$$
\lambda_{1g_{11}}=\cdots = \lambda_{1g_{1n_1}}\quad\text{and}\quad \lambda_{2g_{21}}=\cdots = \lambda_{2g_{2n_2}}.
$$
\end{remark}

\begin{proposition}\label{variante del corolario 3.12} Let $G$, $V$, $f$, $A$, $\alpha$, $V_1$, $V_2$, $g_{11},\dots, g_{1n_1}, g_{21},\dots, g_{2n_2}$, $\hat{\varsigma}$, $\chi_{\varsigma}$, $\hat{\delta}_1$, $\hat{\delta}_2$, $x_1$, $x_2$, $\nu_1$, $\nu_2$, $\lambda_{1g}$ and $\lambda_{2g}$, where $g\in G$, be as in the discussion above Corollary~\ref{coro: second case con P_i en S(V_i)}. Assume that
\begin{xalignat*}{2}
&\lambda_{2g_{11}}=\cdots = \lambda_{2g_{1n_1}}, &&\lambda_{1g_{21}}=\cdots = \lambda_{1g_{2n_2}},\\
&\lambda_{1g_{11}}=\cdots = \lambda_{1g_{1n_1}}, &&\lambda_{2g_{21}}=\cdots =\lambda_{2g_{2n_2}},
\end{xalignat*}
and that conditions~a), b), c) and~d) above that corollary are fulfilled. If
$$
\lambda_{1g_{11}}\lambda_{1g_{21}}=q,\quad \lambda_{2g_{11}}\lambda_{2g_{21}}=q^{-1}\quad \text{and}\quad {}^{g_{ih}}\!x_j=\lambda_{jg_{ih}}x_j,
$$
for $1\le i,j\le 2$ and $1\le h\le n_i$, then:
\begin{enumerate}

\smallskip

\item $\delta_1^l = \delta_2^l = 0$, whenever $q\ne 1$ and $q^l=1$.

\smallskip

\item If $q=1$ or $q$ is not a root of unity, then $P^{(1)}_{g_{1j}}\in \ker \delta_2$ and $P^{(2)}_{g_{2h}}\in \ker \delta_1$ if and only if $P^{(1)}_{g_{1j}}, P^{(2)}_{g_{2h}} \in S(V_1\cap V_2)$.

\smallskip

\item If $q\ne 1$ is a primitive $l$-root of unity, then $P^{(1)}_{g_{1j}}\in \ker \delta_2$ and $P^{(2)}_{g_{2h}}\in \ker \delta_1$ if and only if $P^{(1)}_{g_{1j}}\in S\bigl(k\, x_2^l \oplus (V_1\cap V_2)\bigr)$ and $P^{(2)}_{g_{2h}}\in S\bigl(k\, x_1^l \oplus (V_1\cap V_2)\bigr)$.

\end{enumerate}

\end{proposition}

\begin{proof} Let $\mathbf{x}^{\mathbf{r}} = x_1^{r_1}\cdots x_n^{r_n}$. Using the hypothesis it is easy to check by induction on $s$ that
\begin{align*}
& \delta_1^s(\mathbf{x}^{\mathbf{r}} w_g) = \begin{cases} \displaystyle{\sum_{\mathbf{h}\in \mathds{I}_{n_{\!1}}^s} c_{\mathbf{h}}c'_{\mathbf{h}}} \alpha^s(x_1^{r_1-s} x_2^{r_2}\cdots x_n^{r_n}) w_{g_{1h_s}g_{1h_{s-1}}\cdots g_{1h_1}g} & \text{for $s\le r_1$,}\\ 0 &\text{otherwise,}\end{cases}
\intertext{and}
& \delta_2^s(\mathbf{x}^{\mathbf{r}}w_g) = \begin{cases} \displaystyle{\sum_{\mathbf{h}\in \mathds{I}_{n_{\!2}}^s}} d_{\mathbf{h}}d'_{\mathbf{h}} x_2^{r_2-s} {}^{g_{21}^s}\! \bigr(x_1^{r_1}\! x_3^{r_3}\cdots x_n^{r_n} \bigl) w_{g_{2h_s} g_{2h_{s-1}} \cdots g_{2h_1}g} & \text{for $s\le r_2$,}\\ 0 & \text{otherwise,} \end{cases}
\end{align*}
where
\begin{align*}
& \text{$\mathds{I}_{n_i}^s = \underbrace{\mathds{I}_{n_i}\times\cdots \times \mathds{I}_{n_i}}_{s\text{ times}}$, with $\mathds{I}_{n_i} = \{1,\dots,n_i\}$,}\\[2pt]
& \text{$\alpha^s$ denotes the $s$-fold composition of $\alpha$,}\\
& c_{\mathbf{h}}=\chi_{\varsigma}^s(g)\prod_{k=1}^{s-1}\chi_{\varsigma}^{s-k}(g_{1h_k}) \prod_{k=2}^s \chi_{\alpha}^{k-1}(g_{1h_k}),\\
& c'_{\mathbf{h}} =\Biggl(\prod_{k=0}^{s-1}(r_1-k)_q\Biggr)\Biggl(\prod_{k=1}^s f(g_{1h_k}, g_{1h_{k-1}}\cdots g_{1h_1}g)\Biggr)\Biggl(\prod_{k=1}^s \alpha^{s-1}(P^{(1)}_{g_{1h_k}}) \Biggr),\\
& d_{\mathbf{h}}=\lambda_{2g_{21}}^{sr_2-s(s+1)/2},\\
& d'_{\mathbf{h}}=\Biggl(\prod_{k=0}^{s-1} (r_2-k)_q\Biggr) \Biggl(\prod_{k=1}^s f(g_{2h_k}, g_{2h_{k-1}} \cdots g_{2h_1}g)\Biggr) \Biggl(\prod_{k=0}^{s-1} {}^{g_{21}^k}\! P^{(2)}_{ g_{2h_{s-k}}}\Biggr).
\end{align*}
The result follows easily from these formulas.
\end{proof}

\begin{example} Let $D_u$ be the Dihedral group $D_u := \langle s,\,t \mid s^2,\,t^u,\, stst\rangle$. Then $D_u$ acts on $k[X_1,X_2]$ via
$$
{}^s\! X_1 = -X_1,\quad {}^s\! X_2 = -X_2,\quad {}^t\! X_1 = X_1\quad\text{and}\quad {}^t\! X_2 = X_2.
$$
Let $A = k[X_1,X_2]\#D_u$. We have

\smallskip

\begin{itemize}

\item[-] Assume $u$ is even. Then, there is an $H_1$-module algebra structure on $A$, such that
    \begin{alignat*}{4}
    & \!\qquad\qquad\sigma\!\cdot\! X_1 = X_1,&&\!\quad \sigma\!\cdot\! X_2 = X_2,&&\!\quad \sigma\!\cdot\! w_{t^i} = w_{t^i},&&\!\quad \sigma\!\cdot\! w_{t^i\!s} = - w_{t^i\!s},\\
    & \!\qquad\qquad D_1\!\cdot\! X_1 = w_t+w_{t^{-1}},&&\!\quad D_1\!\cdot\! X_2 = 0,&&\!\quad D_1\!\cdot\! w_{t^i} = 0, &&\!\quad D_1\!\cdot\! w_{t^i\!s} = 0,\\
    & \!\qquad\qquad D_2\!\cdot\! X_1 = 0,&&\!\quad D_2\!\cdot\! X_2 = w_{t^{u/2}},&&\!\quad D_2\!\cdot\! w_{t^i} = 0, &&\!\quad D_2\!\cdot\! w_{t^i\!s} = 0.
    \end{alignat*}

\smallskip

\item[-] There is an $H_{-1}$-module algebra structure on $A$, such that
    \begin{alignat*}{4}
    & \!\qquad\qquad\sigma\!\cdot\! X_1 = X_1,&&\!\quad \sigma\!\cdot\! X_2 = -X_2,&&\!\quad \sigma\!\cdot\! w_{t^i} = w_{t^i},&&\!\quad \sigma\!\cdot\! w_{t^i\!s} = - w_{t^i\!s},\\
    & \!\qquad\qquad D_1\!\cdot\! X_1 = \sum_{i=0}^{u-1} w_{t^i\!s},&&\!\quad D_1\!\cdot\! X_2 = 0,&&\!\quad D_1\!\cdot\! w_{t^i} = 0, &&\!\quad D_1\!\cdot\! w_{t^i\!s} = 0,\\
    & \!\qquad\qquad D_2\!\cdot\! X_1 = 0,&&\!\quad D_2\!\cdot\! X_2 = w_t+w_{t^{-1}},&&\!\quad D_2\!\cdot\! w_{t^i} = 0, &&\!\quad D_2\!\cdot\! w_{t^i\!s} = 0.
    \end{alignat*}

\smallskip

\item[-] Assume $u$ is even. Let $\alpha\colon A\to A$ be the $k$-algebra map defined by $$
    \alpha(Qw_{t^i}) := Qw_{t^i}\quad\text{and}\quad \alpha(Qw_{t^i\!s}) := - Qw_{t^i\!s},
    $$
    and let $s\colon H_1\ot A\to A\ot H_1$ be the transposition associated with $\alpha$. There is an $H_1$-module algebra structure on $A$, such that
    \begin{alignat*}{4}
    & \!\qquad\qquad\sigma\!\cdot\! X_1 = X_1,&&\!\quad \sigma\!\cdot\! X_2 = X_2,&&\!\quad \sigma\!\cdot\! w_{t^i} = w_{t^i},&&\!\quad \sigma\!\cdot\! w_{t^i\!s} = w_{t^i\!s},\\
    & \!\qquad\qquad D_1\!\cdot\! X_1 = w_t+w_{t^{-1}},&&\!\quad D_1\!\cdot\! X_2 = 0,&&\!\quad D_1\!\cdot\! w_{t^i} = 0, &&\!\quad D_1\!\cdot\! w_{t^i\!s} = 0,\\
    & \!\qquad\qquad D_2\!\cdot\! X_1 = 0,&&\!\quad D_2\!\cdot\! X_2 = w_{t^{u/2}},&&\!\quad D_2\!\cdot\! w_{t^i} = 0, &&\!\quad D_2\!\cdot\! w_{t^i\!s} = 0.
    \end{alignat*}

\smallskip

\end{itemize}

\end{example}

\section{Non triviality of the deformations}
%

%
Let $A = S(V)\#_f G$ be as in Section~\ref{Hq-module algebra structures on crossed products}.  By Theorem~\ref{uso de UDF para construir deformaciones} we know that each $H_q$-module algebra $(A,s)$, with $s$ a good transposition, produces to a formal deformation $A_F$ of $A$, which is constructed using the UDF $F=\exp_q(tD_1\!\ot\! D_2)$. The aim of this section is to prove that if $(A,s)$ satisfies the conditions required in Corollary~\ref{coro: second case con P_i en S(V_i)} and $P^{(1)}_{g_{1j}},P^{(2)}_{g_{2h}}\in S(V_1\cap V_2)$ for $1\le j\le n_1$ and $1\le h\le n_2$, then $A_F$ is non trivial. We will prove this showing that its infinitesimal
$$
\Phi(a\ot b) = \delta_1\bigl(\alpha^{-1}(a)\bigr)\delta_2(b),
$$
is not a coboundary. For this we use a complex $\ov{X}^*(A)$, giving the Hochschild cohomology of $A$, which is simpler than the canonical one.

\subsection{A simple resolution}
Given a symmetric $k$-algebra $S:=S(V)$, we consider the differential graded algebra $(Y_*,\partial_*)$ generated by elements $y_v$ and $z_v$, of zero degree, and $\ov{v}$, of degree one, where $v\in V$, subject to the relations
\begin{alignat*}{3}
&z_{\lambda v+w}=\lambda z_v+z_w,\qquad && y_{\lambda v+w}=\lambda y_v+ y_w,\qquad && \ov{v+w}=\lambda \ov{v}+ \ov{w},\\
& y_vy_w = y_wy_v,\qquad && y_vz_w = z_w y_v,\qquad && z_vz_w = z_wz_v,\\
&  \ov{v} y_w = y_w\ov{v},\qquad && \ov{v}z_w = z_w\ov{v}, \qquad && \ov{v}^2=0,
\end{alignat*}
where $\lambda\in k$ and $v,w\in V$, and with differential $\partial$ defined by $\partial(\ov{v}) := \rho_v$, where $\rho_v = z_v - y_v$.

\smallskip

Note that $S$ is a subalgebra of $Y_*$ via the embedding that takes $v$ to $y_v$ for all $v\in V$. This produces an structure of left $S$-module on $Y_*$. Similarly we consider $Y_*$ as a right $S$-module via the embedding of $S$ in $Y_*$ that takes $v$ to $z_v$ for all $v\in V$.
\begin{proposition} Let $\wt{\mu}\colon Y_0\to S$ be the algebra map defined by $\wt{\mu}(y_v) = \wt{\mu}(z_v) := v$ for all $v\in V$. The $S$-bimodule complex
\begin{equation}
\xymatrix{S & Y_0 \lto_-{\wt{\mu}} & Y_1 \lto_-{\partial_1} & Y_2 \lto_-{\partial_2} & Y_3 \lto_-{\partial_3} & Y_4 \lto_-{\partial_4} & Y_5 \lto_-{\partial_5} & \dots \lto_-{\partial_6}}\label{eqqq1}
\end{equation}
is contractible as a left $S$-module complex.
\end{proposition}

\begin{proof} Let $\{x_1,\dots,x_n\}$ be a basis of $V$. We will write $y_i$, $z_i$, $\rho_i$ and $\ov{v}_i$ instead of $y_{x_i}$, $z_{x_i}$, $\rho_{x_i}$ and $\ov{v_{x_i}}$, respectively. A contracting homotopy
$$
\varsigma_0\colon S \to Y_0\quad\text{and}\quad \varsigma_{r+1}\colon Y_r\to Y_{r+1} \quad \text{($r\ge 0$),}
$$
of~\eqref{eqqq1} is given by
\begin{align*}
&\varsigma(1) := 1,\\
&\varsigma\bigl(\rho_{i_1}^{m_1}\ov{v}_{i_1}^{\delta_1}\cdots \rho_{i_l}^{m_l} \ov{v}_{i_l}^{\delta_l} \bigr) := \begin{cases} (-1)^s \rho_{i_1}^{m_1}\ov{v}_{i_1}^{\delta_1} \cdots \rho_{i_{l-1}}^{m_{l-1}} \ov{v}_{i_{l-1}}^{\delta_{l-1}}\rho_{i_l}^{m_l-1} \ov{v}_{i_l} & \text{if $\delta_l = 0$,}\\ 0 & \text{if$\delta_l = 1$,}\end{cases}
\end{align*}
where we assume that $i_1<\cdots<i_l$, $\delta_1+\cdots+\delta_l=s$ and $m_l+\delta_l>0$. In fact, a direct computation shows that

\begin{itemize}

\smallskip

\item[-] $\wt{\mu}\xcirc \sigma^{-1}(1) = \wt{\mu}(1) = 1$,

\smallskip

\item[-] $\varsigma\xcirc\wt{\mu}(1) = \varsigma(1) = 1$ and $\partial\xcirc \varsigma(1) = \partial(0) = 0$,

\smallskip

\item[-] If $\mathbf{x} = \mathbf{x}'\rho_{i_l}^{m_l}$, where $m_l>0$ and $\mathbf{x}' = \rho_{i_1}^{m_1}\cdots \rho_{i_{l-1}}^{m_{l-1}}$ with $i_1<\dots<i_l$, then
$$
\quad\varsigma\xcirc \wt{\mu}(\mathbf{x}) = \varsigma(0) = 0\quad\text{and}\quad\partial \xcirc \varsigma(\mathbf{x}) = \partial(\mathbf{x}'\rho_{i_l}^{m_l-1}\ov{v}_{i_l}) = \mathbf{x},
$$

\smallskip

\item[-] Let $\mathbf{x} = \mathbf{x}'\rho_{i_l}^{m_l} \ov{v}_{i_l}^{\delta_l}$, where $m_l+\delta_l>0$ and $\mathbf{x}' = \rho_{i_1}^{m_1}\ov{v}_{i_1}^{\delta_1}\cdots \rho_{i_{l-1}}^{m_{l-1}} \ov{v}_{i_{l-1}}^{\delta_{l-1}}$ with $i_1<\dots<i_l$ and $\delta_1+\cdots+\delta_l = s>0$. If $\delta_l = 0$, then
\begin{align*}
\quad\qquad &\varsigma\xcirc\partial(\mathbf{x}) = \varsigma\bigl(\partial(\mathbf{x}') \rho_{i_l}^{m_l}\bigr) = (-1)^{s-1}\partial(\mathbf{x}')\rho_{i_l}^{m_l-1} \ov{v}_{i_l},\\
& \partial\xcirc\varsigma(\mathbf{x}) = \partial\bigl((-1)^s\mathbf{x}'\rho_{i_l}^{m_l-1} \ov{v}_{i_l}\bigr) = (-1)^s \partial(\mathbf{x}')\rho_{i_l}^{m_l-1} \ov{v}_{i_l} +\mathbf{x},\\
\intertext{and if $\delta_l = 1$, then}
& \varsigma\xcirc\partial(\mathbf{x}) = \varsigma\bigl(\partial(\mathbf{x}') \rho_{i_l}^{m_l} \ov{v}_{i_l} + (-1)^{s-1} \mathbf{x}'\rho_{i_l}^{m_l+1}\bigr) = \mathbf{x},\\
& \partial\xcirc \varsigma(\mathbf{x}) = \partial(0) = 0.
\end{align*}
\end{itemize}
The result follows immediately from all these facts.
\end{proof}



Let $G$ be a group acting on $V$. We consider $S$ as a $k[G]$-module algebra via the action induced by the one of $G$ on $V$. Let $f\colon k[G]\times k[G]\to k^{\times}$ be a normal cocycle and let $A = S\#_f k[G]$ be the associated crossed product. In the sequel we will use the following

\begin{notation} We let $\ov{k[G]}$ denote $k[G]/k$. Moreover:

\begin{itemize}

\smallskip

\item[-] Given $g_1,\dots,g_s\in \ov{k[G]}$ and $1\le i<j\le s$, we set $\mathbf{g}_{ij} := g_i\ot\cdots\ot g_j$.

\smallskip

\item[-] Given $v_1,\dots,v_r\in V$ and $1\le i<j\le r$, we set $\ov{\mathbf{v}}_{ij} := \ov{v_i}\cdots\ov{v_j}$.

\end{itemize}

\end{notation}

For all $r,s\ge 0$, let
$$
Z_s = (A\ot\ov{k[G]}^{\ot s}) \ot_S A\quad\text{and}\quad X_{rs} = (A\ot\ov{k[G]}^{\ot s})\ot_S Y_r\ot_S A,
$$
where we consider $A\ot\ov{k[G]}^{\ot s}$ as a right $S$-module via
$$
(a_0w_{g_0}\ot \mathbf{g}_{1s})\!\cdot\!a = a_0\, {}^{g_{\!0}\cdots g_{\!s}}\!aw_{g_0}\ot \mathbf{g}_{1s}.
$$
The $X_{rs}$'s and the $Z_s$'s are $A$-bimodules in a canonical way. Note that
$$
Z_s \simeq A\ot\ov{k[G]}^{\ot s}\ot k[G]\quad\text{and}\quad X_{rs} \simeq A\ot\ov{k[G]}^{\ot s}\ot \Lambda^r V\ot A.
$$
In  particular, $X_{rs}$ is a free $A$-bimodule. Consider the diagram of $A$-bimodules and $A$-bimodule maps
$$
\xymatrix{
\vdots \dto^-{-\delta_2} \\
Z_2\dto^-{-\delta_2} & X_{02}\lto_-{\mu_2} & X_{12}\lto_-{d^0_{12}} & \dots\lto_-{d^0_{22}}\\
Z_1\dto^-{-\delta_1} & X_{01}\lto_-{\mu_1} & X_{11}\lto_-{d^0_{11}} & \dots\lto_-{d^0_{21}}\\
Z_0 & X_{00}\lto_-{\mu_0} & X_{10} \lto_-{d^0_{10}} & \dots \lto_-{d^0_{20}},
}
$$
where

\begin{itemize}

\smallskip

\item[-] each $\delta_s$ is defined by
\begin{align*}
\qquad\qquad\delta(1\ot\mathbf{g}_{1s}\ot_S 1) & := w_{g_1}\ot\mathbf{g}_{2s}\ot_S 1\\
& + \sum_{i+1}^{s-1} (-1)^i f(g_i,g_{i+1})\ot\mathbf{g}_{1,i-1}\ot g_ig_{i+1}\ot \mathbf{g}_{i+2,s} \ot_S 1\\
& + (-1)^s 1\ot\mathbf{g}_{1,s-1}\ot_S w_{g_s},
\end{align*}

\smallskip

\item[-] for each $s\ge 0$, the complex $(X_{*s},d_{*s})$ is $(-1)^s$ times $(Y_*,\partial_*)$, tensored over $S$, on the right with $A$ and on the left with $A\ot \ov{k[G]}^{\ot s}$,

\smallskip

\item[-] for each $s\ge 0$, the map $\mu_s$ is defined by
$$
\quad\quad \mu(1\ot\mathbf{g}_{1s}\ot 1) := 1\ot \mathbf{g}_{1s} \ot_S 1.
$$

\smallskip

\end{itemize}
Each row in this diagram is contractible as a left $A$-module. A contracting homotopy
$$
\varsigma^0_{0s}\colon Z_s\to X_{0s}\quad\text{and}\quad \varsigma^0_{r+1,s}\colon X_{rs}\to X_{r+1,s}\quad\text{($r\ge 0$),}
$$
is given by
\begin{align*}
& \varsigma^0(1\ot \mathbf{g}_{1s}\ot_S 1) := 1\ot \mathbf{g}_{1s}\ot 1,\\
& \varsigma^0(1\ot \mathbf{g}_{1s}\ot_S \mathbf{P}\ot_S 1) := (-1)^s 1\ot \mathbf{g}_{1s} \ot_S \varsigma(\mathbf{P})\ot_S 1.
\end{align*}

\medskip

For $r\ge 0$ and $1\le l\le s$, we define $A$-bimodule maps $d^l_{rs}\colon X_{rs}\to X_{r+l-1,s-l}$, recursively on $l$ and $r$, by:
$$
d^l(\mathbf{x}) := \begin{cases}
\varsigma^0\xcirc\delta\xcirc\mu(\mathbf{x}) &\text{if $l=1$ and $r=0$,}\\
- \varsigma^0\xcirc d^1\xcirc d^0(\mathbf{x}) &\text{if $l=1$ and $r>0$,}\\
- \sum_{j=1}^{l-1}\varsigma^0\xcirc d^{l-j}\xcirc d^j(\mathbf{x})&\text{if $1<l$ and $r=0$,}\\
- \sum_{j=0}^{l-1}\varsigma^0\xcirc d^{l-j}\xcirc d^j(\mathbf{x}) &\text{if $1<l$ and $r>0$,}
\end{cases}
$$
for $\mathbf{x} = 1 \ot\mathbf{g}_{1s}\ot \ov{\mathbf{v}}_{1r}\ot 1$.

\begin{theorem}\label{res nuestra} There is a resolution of $A$ as an $A$-bimodule
\begin{equation*}
\xymatrix{A & X_0\lto_{-\mu} & X_1  \lto_{d_1}  &X_2 \lto_{d_2}  & X_3 \lto_{d_3}  & X_4 \lto_{d_4} & \lto_{d_5}\dots,}
\end{equation*}
where $\mu\colon X_{00}\to A$ is the multiplication map,
$$
X_n = \bigoplus_{r+s=n} X_{rs}\quad\text{and}\quad d_n = \sum^n_{l=1} d^l_{0n} + \sum_{r=1}^n \sum^{n-r}_{l=0} d^l_{r,n-r}.
$$
\end{theorem}

\begin{proof} See~\cite[Appendix A]{G-G2}.
\end{proof}

\begin{proposition}\label{Formulas para d^l} The maps $d^l$ vanish for all $l\ge 2$. Moreover
\begin{align*}
d^1(1\ot\mathbf{g}_{1s}\ot \ov{\mathbf{v}}_{1r}\ot 1) &= w_{g_1}\ot\mathbf{g}_{2s} \ot \ov{\mathbf{v}}_{1r}\ot 1\\
& + \sum_{i=1}^{s-1} (-1)^if(g_i,g_{i+1})\ot\mathbf{g}_{1,i-1}\ot g_ig_{i+1}\ot
\mathbf{g}_{i+2,s} \ot \ov{\mathbf{v}}_{1r}\ot 1\\
& + (-1)^s 1\ot \mathbf{g}_{1,s-1}\ot\ov{{}^{g_{\!s}}\!v_1}\cdots\ov{{}^{g_{\!s}}\!v_r} \ot w_{g_s}.
\end{align*}
In particular, $(X_*,d_*)$ is the total complex of the double complex
$$
\xymatrix{
\vdots \dto^-{d^1_{03}} & \vdots \dto^-{d^1_{13}} & \vdots \dto^-{d^1_{23}}\\
X_{02} \dto^-{d^1_{02}} & X_{12} \lto_-{d^0_{12}}\dto^-{d^1_{12}} & X_{22} \lto_-{d^0_{22}} \dto^-{d^1_{22}} & \dots \lto_-{d^0_{32}}\\
X_{01} \dto^-{d^1_{01}} & X_{11} \lto_-{d^0_{11}}\dto^-{d^1_{11}} & X_{21} \lto_-{d^0_{21}} \dto^-{d^1_{21}} & \dots \lto_-{d^0_{31}}\\
X_{00} & X_{10} \lto_-{d^0_{10}} & X_{20} \lto_-{d^0_{20}} & \dots \lto_-{d^0_{30}},
 }
$$
\end{proposition}

\begin{proof} The computation of $d_{rs}^1$ can be obtained easily by induction on $r$, using that
\begin{align*}
& d^1(\mathbf{x}) = \varsigma^0\xcirc\delta\xcirc\mu(\mathbf{x})\quad\text{for $\mathbf{x} = 1\ot \mathbf{g}_{1s}\ot 1$,}
\intertext{and}
& d^1(\mathbf{x}) = - \varsigma^0\xcirc d^1\xcirc d^0(\mathbf{x})\quad\text{for $r\ge 1$ and $\mathbf{x} = 1 \ot\mathbf{g}_{1s}\ot \ov{\mathbf{v}}_{1r}\ot 1$.}
\end{align*}
The assertion for $d_{rs}^l$, with $l\ge 2$, follows by induction on $l$ and $r$, using the recursive definition of $d_{rs}^l$.
\end{proof}

\subsection{A comparison map}
Let $\ov{A} = A/k$. In this subsection we introduce and study a comparison map from $(X_*,d_*)$ to the canonical normalized Hochschild resolution $(A\ot\ov{A}^*\ot A,b'_*)$. It is well known that there is an $A$-bimodule homotopy equivalence
$$
\theta_*\colon (X_*,d_*)\to (A\ot\ov{A}^*\ot A,b'_*)
$$
such that $\theta_0 =\ide_{A\ot A}$. It can be recursively defined by $\theta_0 := \ide_{A\ot A}$ and
$$
\theta(\mathbf{x}) := (-1)^{r+s}\theta\xcirc d(\mathbf{x})\ot 1\quad \text{for $\mathbf{x} = 1\ot \mathbf{g}_{1s}\ot \ov{\mathbf{v}}_{1r}\ot 1$ with $r+s\ge 1$.}
$$
Next we give a closed formula for $\theta_*$. In order to establish this result we need to introduce a new notation. We recursively define $(w_{g_1}\ot \cdots\ot w_{g_s}) * (P_1\ot\cdots\ot P_r)$ by

\begin{itemize}

\smallskip

\item[-] $(w_{g_1}\ot \cdots\ot w_{g_s}) * (Q_1\ot\cdots\ot Q_r) := (Q_1\ot\cdots\ot Q_r)$ if $s=0$,

\item[-] $(w_{g_1}\ot \cdots\ot w_{g_s}) * (Q_1\ot\cdots\ot Q_r) := (w_{g_1}\ot \cdots\ot w_{g_s})$ if $r=0$,

\smallskip

\item[-] If $r,s\ge 1$, then $(w_{g_1}\ot \cdots\ot w_{g_s}) * (Q_1\ot\cdots\ot Q_r)$ equals
$$
\qquad\qquad\sum_{i=0}^r (-1)^i (w_{g_1}\ot\cdots\ot w_{g_{s-1}})*({}^{g_{\!s}}\! Q_1\ot \cdots \ot {}^{g_{\!s}}\! Q_i)\ot w_{g_s}\ot Q_{i+1}\ot\cdots\ot Q_r.
$$

\smallskip

\end{itemize}

\begin{proposition}\label{prop 3.1} We have:
$$
\theta\bigl(1\ot\mathbf{g}_{1s}\ot \ov{\mathbf{v}}_{1r}\ot 1\bigr) = (-1)^r\sum_{\tau\in \mathfrak{S}_r}\sg(\tau)\ot\bigl(w_{g_1}\ot\cdots\ot w_{g_s}\bigr)*\mathbf{v}_{\tau(1r)}\ot 1,
$$
where $\mathfrak{S}_r$ is the symmetric group in $r$ elements and $\mathbf{v}_{\tau(1r)} = v_{\tau(1)}\ot\cdots\ot v_{\tau(r)}$.
\end{proposition}

\begin{proof} We proceed by induction on $n=r+s$. The case $n=0$ is obvious. Suppose that $r+s =n$ and the result is valid for $\theta_{n-1}$. By the recursive definition of $\theta$ and Theorem~\ref{res nuestra},
\begin{align*}
\theta(1\ot \mathbf{g}_{1s}\ot \ov{\mathbf{v}}_{1r}\ot 1) & = (-1)^n\theta\xcirc d(1\ot \mathbf{g}_{1s} \ot \ov{\mathbf{v}}_{1r}\ot 1)\ot 1\\
& = (-1)^n\theta\xcirc (d^0+d^1)(1\ot\mathbf{g}_{1s}\ot\ov{\mathbf{v}}_{1r}\ot 1)\ot 1\\
&= \sum_{i=1}^r (-1)^{i+r}\theta({}^{g_{\!1}\cdots g_{\!s}}\! v_i\ot\mathbf{g}_{1s}\ot \ov{\mathbf{v}}_{1,i-1} \ov{\mathbf{v}}_{i+1,r} \ot 1)\ot 1\\
& - \sum_{i=1}^r(-1)^{i+r}\theta(1\ot\mathbf{g}_{1s}\ot \ov{\mathbf{v}}_{1,i-1} \ov{\mathbf{v}}_{i+1,r} \ot v_i)\ot 1\\
& + (-1)^n\theta(w_{g_1}\ot\mathbf{g}_{2s}\ot\ov{\mathbf{v}}_{1r}\ot 1)\ot 1\\
& + \sum_{i=1}^{s-1} (-1)^{n+i}\theta(1\ot\mathbf{g}_{1,i-1}\ot g_ig_{i+1}\ot\mathbf{g}_{i+1,s} \ot \ov{\mathbf{v}}_{1r} \ot 1) \ot 1\\
& +(-1)^r\theta(1\ot\mathbf{g}_{1,s-1}\ot \mathbf{g}_{1,s-1}\ot\ov{{}^{g_{\!s}}\!v_1}\cdots \ov{{}^{g_{\!s}}\!v_r}\ot w_{g_s})\ot 1.
\end{align*}
The desired result follows now from the inductive hypothesis.
\end{proof}

\subsection{The Hochschild cohomology}
Let $M$ be an $A$-bimodule and $A^e$ the enveloping algebra of $A$. Applying the functor $\Hom_{A^e}(-,M)$ to $(X_{**},d^0_{**},d^1_{**})$ and using the identifications
$$
\Hom_{A^e}(X_{rs},M) \simeq \Hom_k(\ov{k[G]}^{\ot s}\ot \Lambda^r V,M)
$$
we obtain the double complex
$$
\xymatrix{
\vdots  & \vdots  & \vdots \\
\ov{X}^{02}\uto_-{\ov{d}_1^{03}} \rto^-{\ov{d}_0^{12}} & \ov{X}^{12} \uto_-{\ov{d}_1^{13}} \rto^-{\ov{d}_0^{22}}& \ov{X}^{22} \rto^-{\ov{d}_0^{32}} \uto_-{\ov{d}_1^{23}} & \dots \\
\ov{X}^{01}\uto_-{\ov{d}_1^{02}} \rto^-{\ov{d}_0^{11}} & \ov{X}^{11} \uto_-{\ov{d}_1^{12}} \rto^-{\ov{d}_0^{21}}& \ov{X}^{21} \rto^-{\ov{d}_0^{31}}\uto_-{\ov{d}_1^{22}} & \dots \\
\ov{X}^{00} \uto_-{\ov{d}_1^{01}}\rto^-{\ov{d}_0^{10}}& \ov{X}^{10} \uto_-{\ov{d}_1^{11}} \rto^-{\ov{d}_0^{20}} & \ov{X}^{20} \uto_-{\ov{d}_1^{21}}\rto^-{\ov{d}_0^{30}} & \dots,
}
$$
where
\allowdisplaybreaks
\begin{align*}
& \ov{X}^{rs} = \Hom_k(\ov{k[G]}^{\ot s}\ot \Lambda^r V,M),\\
& \ov{d}_0(\varphi)(\mathbf{g}_{1s}\ot\ov{\mathbf{\mathbf{v}}}_{1,r+1})=\sum_{i=1}^{r+1} (-1)^{s+i+1} \varphi(\mathbf{g}_{1s}\ot \ov{\mathbf{\mathbf{v}}}_{1,i-1} \ov{\mathbf{\mathbf{v}}}_{i+1,r+1})v_i\\
&\phantom{\ov{d}_0(\varphi)(\mathbf{g}_{1s}\ot\ov{\mathbf{\mathbf{v}}}_{1,r+1})} + \sum_{i=1}^{r+1} (-1)^{s+i}\, {}^{g_{\!1}\cdots g_{\!s}}\!v_i \varphi(\mathbf{g}_{1s}\ot \ov{\mathbf{\mathbf{v}}}_{1,i-1}\ov{\mathbf{\mathbf{v}}}_{i+1,r+1}),\\
& \ov{d}_1(\varphi)(\mathbf{g}_{1,s+1}\ot\ov{\mathbf{\mathbf{v}}}_{1r}) = w_{g_1} \varphi(\mathbf{g}_{2,s+1}\ot\ov{\mathbf{v}}_{1r})\\
&\phantom{\ov{d}_1(\varphi)(\mathbf{g}_{1,s+1}\ot\ov{\mathbf{\mathbf{v}}}_{1r})} + \sum_{i=1}^s (-1)^i f(g_i,g_{i+1})\varphi (\mathbf{g}_{1,i-1}\ot g_ig_{i+1}\ot \mathbf{g}_{i+1,s+1}\ot \ov{\mathbf{v}}_{1r})\\
&\phantom{\ov{d}_1(\varphi)(\mathbf{g}_{1,s+1}\ot\ov{\mathbf{\mathbf{v}}}_{1r})} + (-1)^{s+1} \varphi(\mathbf{g}_{1s}\ot\ov{{}^{g_{\!s+1}}\!v_1}\cdots \ov{{}^{g_{\!s+1}}\!v_r}) w_{g_{s+1}},
\end{align*}
whose total complex $\ov{X}^*(M)$ gives the Hochschild cohomology $\Ho^*(A,M)$ of $A$ with coefficients in $M$. The comparison map $\theta_*$ induces a quasi-isomorphism
$$
\ov{\theta}^*\colon \bigl(\Hom_k(\ov{A}^*,M),\ov{b}^*\bigr) \to \ov{X}^*(M).
$$
It is immediate that
$$
\ov{\theta}(\varphi)(\mathbf{g}_{1s}\ot \ov{\mathbf{v}}_{1r})= (-1)^r\sum_{\tau\in \mathfrak{S}_r} \sg(\tau)\varphi \bigl((w_{g_1}\ot\cdots\ot w_{g_s})* \mathbf{v}_{\tau(1r)} \bigr).
$$
where $\mathfrak{S}_r$ is the symmetric group in $r$ elements and $\mathbf{v}_{\tau(1r)} = v_{\tau(1)}\ot\cdots\ot v_{\tau(r)}$.

\smallskip

From now on we take $M = A$ and we write $\HH^*(A)$ instead of $\Ho^*(A,A)$.

\subsection{Proof of the main result} We are ready to prove that the cocycle $\Phi$ is non trivial. For this it is sufficient to show that $\ov{\theta}(\Phi)$ is not a coboundary. Let $x_1,\dots,x_n$, $P^{(1)}_{g_{11}},\dots,P^{(1)}_{g_{1n_1}}$, $P^{(2)}_{g_{21}},\dots, P^{(2)}_{g_{2n_2}}$, $g_{11},\dots g_{1n_1}$ and $g_{21},\dots g_{2n_2}$ be as in Corollary~\ref{coro: second case con P_i en S(V_i)}. A direct computation, using the formulas for $\delta_1$ and $\delta_2$ obtained in the proof of Proposition~\ref{variante del corolario 3.12}, shows that
$$
\ov{\theta}(\Phi)(g\ot \ov{v}) = 0\quad\text{and}\quad  \ov{\theta}(\Phi)(g\ot h) = 0.
$$
for $g,h\in G$ and $v\in V$, and that
$$
\ov{\theta}(\Phi)(\ov{x_1}\,\ov{x_2}) = \sum_{j=1}^{n_1} \sum_{h=1}^{n_2}\chi_{\alpha}^{-1} (g_{1j}) f(g_{1j},g_{2h}) \alpha^{-1}(P^{(1)}_{g_{1j}}) {}^{g_{1j}}\!P^{(2)}_{g_{2h}} w_{g_{1j} g_{2h}}
$$
and
$$
\ov{\theta}(\Phi)(\ov{x_i}\,\ov{x_j}) = 0\quad\text{for $1\le i<j\le n$ with $(i,j)\ne (1,2)$.}
$$
We next prove that $\ov{\theta}(\Phi)$ is not a coboundary. Let $\varphi_0\in \ov{X}_{01}$ and $\varphi_1\in \ov{X}_{10}$. By definition
\begin{align*}
& \ov{d}_1(\varphi_0)(g\ot h) = w_g\varphi_0(h) - f(g,h)\varphi_0(gh) + \varphi_0(g)w_h,\\
& \ov{d}_0(\varphi_0)(g\ot \ov{v}) = {}^g\!v\varphi_0(g) - \varphi_0(g)v,\\
& \ov{d}_1(\varphi_1)(g\ot \ov{v}) = w_g\varphi_1(\ov{v}) - \varphi_1(\ov{{}^g\!v})w_g,\\
& \ov{d}_0(\varphi_1)(\ov{v_1}\,\ov{v_2}) = \varphi_1(\ov{v_2})v_1 - v_1\varphi_1(\ov{v_2}) + v_2\varphi_1(\ov{v_1}) - \varphi_1(\ov{v_1})v_2,
\end{align*}
and so $\ov{\theta}(\Phi)$ is a coboundary if and only if there exist $\varphi_0$ and $\varphi_1$ such that
\begin{align*}
& w_g\varphi_0(h) - f(g,h)\varphi_0(gh) + \varphi_0(g)w_h = 0 && \text{for all $g,h\in G$,}\\
& {}^g\!v\varphi_0(g) - \varphi_0(g)v + w_g\varphi_1(\ov{v}) - \varphi_1(\ov{{}^g\!v})w_g = 0 && \text{for all $g\in G$ and $v\in V$,}\\
& [\varphi_1(\ov{x_j}),x_i]  + [x_j,\varphi_1(\ov{x_i})] = 0 && \text{for all $i<j$ with $(i,j)\ne (1,2)$,}
\end{align*}
where, as usual, $[a,b] = ab - ba$, and
$$
[\varphi_1(\ov{x_2}),x_1] + [x_2,\varphi_1(\ov{x_1})] = \sum_{j=1}^{n_1} \sum_{h=1}^{n_2} \chi_{\alpha}^{-1} (g_{1j}) f(g_{1j},g_{2h}) \alpha^{-1}(P^{(1)}_{g_{1j}}) {}^{g_{1j}}\! P^{(2)}_{g_{2h}} w_{g_{1j}g_{2h}}.
$$
But, since $w_gx_j = {}^g\!x_j w_g$,
$$
w_{g_{1j}g_{2h}}x_1 = f(g_{1j},g_{2h})^{-1} w_{g_{1j}}w_{g_{2h}}x_1 = qx_1\quad\text{and}\quad w_{g_{1j}g_{2h}}x_2 = q^{-1}x_2,
$$
if
$$
\varphi_1(\ov{x_1}) = \sum_{g\in G} Q^{(1)}_g w_g\quad\text{and}\quad \varphi_1(\ov{x_2}) = \sum_{g\in G} Q^{(2)}_g w_g,
$$
then necessarily
$$
\sum_{g\in \Upsilon}(q-1)\bigl(x_1 Q^{(2)}_g + q^{-1}x_2 Q^{(1)}_g\bigr) w_g= \sum_{j=1}^{n_1} \sum_{h=1}^{n_2} D_{jh} \alpha^{-1}(P^{(1)}_{g_{1j}}) {}^{g_{1j}}\!P^{(2)}_{g_{2h}} w_{g_{1j} g_{2h}},
$$
where
$$
D_{jh}=\chi_{\alpha}^{-1} (g_{1j}) f(g_{1j},g_{2h}) \quad\text{and}\quad\Upsilon = \{g_{1j}g_{2h}:1\le j\le n_1\text{ and } 1\le h\le n_2\},
$$
which is impossible because $\alpha^{-1}(P^{(1)}_{g_{1j}}) {}^{g_{1j}}\!P^{(2)}_{g_{2h}}\in k[x_3,\dots,x_n] \setminus\{0\}$.


\end{document}